\documentclass[11pt]{article}
\usepackage{amsmath,amsthm,amssymb,amscd}

\usepackage{graphicx}
\usepackage{epstopdf}
\usepackage{color}
\usepackage{url}

\theoremstyle{definition}\newtheorem{definition}{Definition}
\theoremstyle{plain}\newtheorem{theorem}[definition]{Theorem}
\theoremstyle{definition}
\theoremstyle{definition}\newtheorem{proposition}[definition]{Proposition}
\theoremstyle{definition}
\theoremstyle{definition}
\theoremstyle{definition}
\theoremstyle{definition}\newtheorem{lemma}[definition]{Lemma}
\theoremstyle{definition}

\newcommand{\cR}{{\mathcal R}}

\newcommand{\cN}{{\mathcal N}}
\newcommand{\R}{\mathbb{R}}
\newcommand{\N}{\mathbb{N}}

\newcommand{\xad}{x_\alpha^\delta}
\newcommand{\xd}{x^\dag}
\newcommand{\xa}{x_\alpha}

\newcommand{\xiad}{\xi_\alpha^\delta}
\newcommand{\xia}{\xi_\alpha}

\newcommand{\asa}{A^\ast A}
\newcommand{\basa}{(A^\ast A)}
\newcommand{\bigo}{\mathcal{O}}

\begin{document}

\title{A new interpretation of (Tikhonov) regularization}

\author{Daniel Gerth\footnote{Technische Universit\"at Chemnitz, Faculty for Mathematics, D-09107 Chemnitz, Germany, daniel.gerth@mathematik.tu-chemnitz.de}}


%
%

\maketitle

\begin{abstract}
Tikhonov regularization with square-norm penalty for linear forward operators has been studied extensively in the literature. However, the results on convergence theory are based on technical proofs and difficult to interpret. It is also often not clear how those results translate into the discrete, numerical setting. In this paper we present a new strategy to study the properties of a regularization method on the example of Tikhonov regularization. The technique is based on the observation that Tikhonov regularization approximates the unknown exact solution in the range of the adjoint of the forward operator. This is closely related to the concept of approximate source conditions, which we generalize to describe not only the approximation of the unknown solution, but also noise-free and noisy data; all from the same source space. Combining these three approximation results we derive the well-known convergence results in a concise way and improve the understanding by tightening the relation between concepts such as convergence rates, parameter choice, and saturation. The new technique is not limited to Tikhonov regularization, it can be applied also to iterative regularization, which we demonstrate by relating Tikhonov regularization and Landweber iteration. Because the Tikhonov functional is no longer the centrepiece of the analysis, we can show that Tikhonov regularization can be used for oversmoothing regularization. All results are accompanied by numerical examples.
\end{abstract}

\section{Introduction}
The regularization of linear ill-posed problems, in particular via classical Tikhonov regularization with square-norm penalty, has been discussed extensively in the literature and is considered well-understood. 
However, since this understanding is often based on technical proofs and conditions, an intuitive interpretation is often missing. The goal of this paper is to introduce such an interpretation that can be used as a new basis for regularization theory not only for classical Tikhonov regularization but also other methods, even iterative ones such as Landweber iteration as we briefly demonstrate. However, for the sake of brevity, we focus mostly on Tikhonov regularization. Our new approach essentially reduces regularization to a single basic principle, namely, that the exact solution can be approximated arbitrarily well in a space of higher smoothness. This point of view decouples exact data, noise, and the regularization method which opens up a more universal and general approach for the analysis of ill-posed problems.
In contrast to the classical theory, which relies heavily on the Tikhonov functional, we work almost exclusively with the first-order optimality condition since it can be directly related to the approximation properties. Based on that, we recover well-known results on, e.g., convergence rates, parameter choice, and saturation. While our method yields neither improvements for those properties nor computational benefits, it explains known connections between the theoretical results in a much more comprehensible way and opens new connections. An additional advantage is that the treatment of discrete ill-posed problems is naturally included in the approximation approach, which establishes a new and simple way to relate discrete and infinite-dimensional ill-posed problems. 

In order to explain our concept in more detail and to point out some questions on the classical theory, we introduce in the following some notation and the basic setting. Let $X$ be the Hilbert space containing the sought-after solution $x^\dag$, and $Y$ be the data space containing the data $y=Ax^\dag$, where $A:X\rightarrow Y$ is assumed to be a compact linear operator. The compactness implies that $A$ has non-closed range, $\cR(A)\neq\overline{\cR(A)}$. While there are operators with non-closed ranged that are not compact (called strictly singular operators), we confine ourselves here to compactness since it allows to use the singular system of $A$ for computations, and constitutes a natural limit for finite dimensional (and thus necessarily compact) approximations to $A$ used numerically. Let $\{\sigma_i,u_i,v_i\}_{i=1}^\infty$ be the singular system of $A$ where the $\{u_i\}_{i=1}^\infty$ form an orthonormal basis (ONB) for $\overline{\cR(A)}$, $\{v_i\}_{i=1}^\infty$ form an ONB for $\overline{\cR(A^\ast)}$ and the singular values $\{\sigma_i\}$ accumulate at zero. We recall that the relations $Av_i=\sigma_i u_i$ and $A^\ast u_i=\sigma_i v_i$ hold for all $i\in \N$, and that
\[
Ax= \sum_{i=1}^\infty \langle Ax,u_i\rangle u_i=\sum_{i=1}^\infty \sigma_i\langle x,v_i\rangle u_i .
\]

For the remainder of the paper we are concerned with the solution of the ill-posed linear operator equation
\begin{equation}\label{eq:problem}
y=Ax,
\end{equation}
which for notational convenience we will consider to be scaled such that $\|A\|=\sigma_1=1$. In practice we only have access to noisy data $y^\delta$, for which we use the convention
\begin{equation}\label{eq:noise}
\|y-y^\delta\|\leq \delta
\end{equation}
with some $0<\delta\leq\delta_0$. To stabilize the recovery we consider Tikhonov regularization in its classical formulation, i.e., the approximate solution to (\ref{eq:problem}) is obtained as
\begin{equation}\label{eq:tikh_classic}
x_\alpha^\delta=\mathrm{argmin}_{x\in X} \left\lbrace \frac{1}{2}\|Ax-y^\delta\|^2+\frac{\alpha}{2}\|x\|^2\right\rbrace.
\end{equation}
In the case of noise-free data we denote the regularized solutions by $\xa$, i.e.,
\begin{equation*}\label{eq:tikh_classic_nonoise}
x_\alpha=\mathrm{argmin}_{x\in X} \left\lbrace \frac{1}{2}\|Ax-y\|^2+\frac{\alpha}{2}\|x\|^2\right\rbrace.
\end{equation*}
We recall that the first-order optimality condition for (\ref{eq:tikh_classic}) is
\begin{equation}\label{eq:tikh_optcond}
A^\ast(A\xad-y^\delta)+\alpha\xad=0,
\end{equation}
and hence the solution to (\ref{eq:tikh_classic}) is given by
\begin{equation}\label{eq:tikh_explicit}
\xad=(\asa+\alpha I)^{-1}A^\ast y^\delta.
\end{equation}
Using the singular system of $A$, we can write
\begin{equation*}\label{eq:tikh_svd}
\xad =\sum_{i=1}^\infty \frac{\sigma_i}{\sigma_i^2+\alpha}\langle y^\delta,u_i\rangle v_i.
\end{equation*}
In practice, the most important question is how to choose the regularization parameter $\alpha$, which, with $A$, $\xd$, and $y^\delta$ fixed, determines the reconstruction error $\|\xad-\xd\|$. This is closely related to the theory of convergence rates, where one looks for an \textit{index function} $\varphi$, i.e., $\varphi:[0,\infty)\rightarrow \R_+$ is continuous and monotonically increasing with $\varphi(0)=0$, such that for a suitable choice of the regularization parameter $\alpha$
\begin{equation*}\label{eq:rates_prototype}
\|\xad-\xd\|\leq \varphi(\delta) \qquad \mbox{for all } 0<\delta\leq \delta_0.
\end{equation*}
Without any imposed restriction, no such $\varphi$ exists. In order to find a convergence rate, a smoothness relation between operator $A$ and solution $x^\dag$ has to be established. The classical assumption in this regard is a \textit{source condition}, postulating the existence of a parameter $\mu>0$ such that 
\begin{equation}\label{eq:sc}
x^\dag\in\cR((A^\ast A)^\mu).
\end{equation}
We prefer in this paper the slightly modified condition
\begin{equation}\label{eq:sc_open}
x^\dag\in \bigcap_{\kappa<\mu} \cR((A^\ast A)^\kappa),
\end{equation}
which in the singular system is implied if
\begin{equation}\label{eq:sc_sum}
\sum_{n=k}^\infty \langle \xd,v_n\rangle^2=\bigo(\sigma_k^{4\mu})
\end{equation}
for $k\rightarrow \infty$, see \cite{Neubauer}. Assuming (\ref{eq:sc_open}) is satisfied with some $0<\mu<1$, one can show (see, e.g., \cite{EHN}) that
\begin{equation}\label{eq:opt_rate1}
\|\xad-\xd\|\leq C \delta^{\frac{2\mu}{2\mu+1}}
\end{equation}
is, even without considering any regularization, the best obtainable convergence rate in the sense that the exponent can not be increased any further. Tikhonov regularization yields this rate if (\ref{eq:sc}) holds with $0<\mu\leq 1$ and if $\alpha$ is chosen appropriately, for example if
\begin{equation}\label{eq:apriori1}
C_1\delta^{\frac{2}{2\mu+1}}\leq\alpha\leq C_2\delta^{\frac{2}{2\mu+1}}
\end{equation}
with appropriate constants $0<C_1\leq C_2< \infty$ or with the discrepancy principle, i.e., selecting $\alpha$ such that
\begin{equation}\label{eq:dp1}
\alpha^\ast=\sup \{\alpha>0: \|A\xad-y^\delta\|\leq \tau \delta, \tau>1\}
\end{equation}
However, the latter only yields (\ref{eq:opt_rate1}) for $0<\mu<\frac{1}{2}$ in (\ref{eq:sc}). 

The results presented above are well-known but rely on technical proofs that confirm their correctness, and are difficult to comprehend and interpret on an intuitive level. What do the exponents in the convergence rate prototype describe, why are they universal in the sense that no regularization method can improve them? We show that it is a result of the difference in the approximation properties of $\xd$ and the data $y$ from the same source space. This includes the statement that for well-posed problems, the optimal convergence rate is of order $\delta$. For Tikhonov regularization, this source space is naturally the domain $\mathrm{dom}(A^\ast)$ of the adjoint of $A$ because we show $\xa,\xad\in \cR(A^\ast)$ unconditionally. Why does the discrepancy principle yield the rates (\ref{eq:opt_rate1}) only for $0<\mu<\frac{1}{2}$? We show that this is due to $\mathrm{dom}(A^\ast)$ becoming non-informative for $\mu>\frac{1}{2}$ since then $\xd\in\cR((A^\ast A)^\mu)\subset A^\ast$ instead of $\xd\in\cR((A^\ast A)^\mu)\supset \cR(A^\ast)$ in the case $0<\mu<\frac{1}{2}$. Using this argument, we can show that the a-priori parameter choice (\ref{eq:apriori1}) and the discrepancy principle (\ref{eq:dp1}) (and further any parameter choice rule yielding order optimal convergence rates) coincide up to a constant, because their only role is to realize a specific growth of the norm of the source element which, in the case of noisy data, corresponds to ensure $\|A\xad-y^\delta\|\approx \delta$. Clearly, the follow-up question must be why the a-priori choice then yields the rate (\ref{eq:opt_rate1}) even for $\frac{1}{2}\leq \mu\leq 1$. We explain this by considering $\xa\in \cR(A^\ast A)$. It is also this representation that restricts higher convergence rates, because it limits the choice of the regularization parameters by forcing $\|\xa-\xd\|/\alpha$ to be bounded from below independent of $\mu$. Therefore, the new approach yields a simple way of analysing saturation of regularization methods, which we use construct a variant of Tikhonov regularization that can yield the rate (\ref{eq:opt_rate1}) for arbitrary fixed $0<\mu<\infty$ using an a-priori parameter choice or for $0<\mu-\frac{1}{2}<\infty$ with the discrepancy principle. The key to achieve this is to force the regularized solution to lie in spaces smoother then $\cR((A^\ast A)^\mu)$.
This forced smoothness relation is also what enables oversmoothing regularization, i.e., Tikhonov regularization under the assumption that $\|\xd\|=\infty$. The difficulty in classical analysis is that in this situation the estimate
\[
\|A\xad-y^\delta\|^2+\alpha\|\xad\|^2\leq \|A\xd-y^\delta\|^2+\alpha\|\xd\|^2
\]
becomes meaningless as the right hand side is infinite. In our new approach, there is no difference between oversmoothing regularization $\|\xd\|=\infty$ and the classical setting $\|\xd\|<\infty$ because the first-order optimality condition is independent of the norm of $\xd$. To demonstrate the applicability of our new approach we consider regularization in Hilbert scales, for which we recover and slightly extend the classical results based on the same theorem and subsequent analysis as used for classical Tikhonov regularization, just with slightly adjusted parameters. 
Finally we briefly show that the principle of approximating $\xd$ in smoother spaces also holds for iterative regularization methods. We consider Landweber method as a particular example and demonstrate that reconstruction errors and residuals obtained with Landweber iteration (for many individual iterates) and Tihkonov regularization (for many parameters $\alpha$) are almost identical when plotted against the norm of the respective source element in $\mathrm{dom}(A^\ast)$.

The paper is structured as follows. In Section \ref{sec:main} we show $\xa,\xad\in\cR(A^\ast)$ unconditionally and illustrate in an example that this immediately translates to effects observable in numerical computations. In Section \ref{sec:appsc} we discuss the concept of approximate source conditions, which had previously been used to analyse Tikhonov regularization. The novelty in this paper is that the approximate source conditions are extended to approximate powers $\basa^\nu\xd$ including the approximation of the data, and to approximate the noise. We show that the optimal convergence rate for any regularization method is obtained by combining approximate source conditions for $\xd$ and $A\xd$. Based on this we derive the well-known results on convergence rates and parameter choice by relating the magnitude of the corresponding source elements. We also show that the infinite dimensional approximation properties are retained to a large extent in the discrete setting. The approximation-based view motivates a discussion of saturation in Section \ref{sec:sat}. Having understood saturation properly, we demonstrate how increasing the smoothness of the approximate solutions by adjusting Tikhonov regularization yields higher, and precisely controllable and predictable saturation properties, see Section \ref{sec:higherTikh}. Having understood that a main principle of regularization is that the approximate solutions be smoother than $\xd$, we discuss the case of oversmoothing regularization in Section \ref{sec:os}. The same principle can also be transferred to iterative regularization method, which we demonstrate on the example of Landweber iteration. In terms of the growth of the source element, there is only a negligible difference between Tikhonov-regularization and Landweber iteration.

\section{Main observation}\label{sec:main} 
All further discussions are based on the following observation.

\begin{proposition}\label{thm:bas}
The Tikhonov-regularized solutions to (\ref{eq:tikh_classic}) can be written in the form
\begin{equation}\label{eq:tikh_sourcerep1}
\xad=A^\ast \frac{A\xad-y^\delta}{\alpha}, \qquad \xa=A^\ast \frac{A\xa-y}{\alpha}
\end{equation}
i.e., $\xad,\xa \in \cR(A^\ast)$. Further we have the representation
\begin{equation}\label{eq:tikh_sourcerep2}
\xad=A^\ast A \frac{\xad-\xd}{\alpha}-A^\ast\frac{y-y^\delta}{\alpha},\qquad \xa=A^\ast A \frac{\xa-\xd}{\alpha}
\end{equation}
which, at least in the noise-free case, implies $\xa\in \cR(A^\ast A)$.
\end{proposition}
\begin{proof}
All equations are reorderings of the first-order optimality condition (\ref{eq:tikh_optcond}) for noisy and noise-free data, respectively.
\end{proof}

The implicit equations (\ref{eq:tikh_sourcerep1}) and (\ref{eq:tikh_sourcerep2}) have no advantage for the calculation of the approximate solutions $\xad,\xa$. For this, the explicit formula (\ref{eq:tikh_explicit}) is much more suited. Instead, Proposition \ref{thm:bas} gives a straight-forward description of the smoothness of the approximate solutions in Tikhonov regularization: it is $\xad,\xa \in \cR(A^\ast)$ unconditionally. Note that this is stronger than the smoothness $\|\xad\|<\infty$, or $\xad\in X$, as implied by the minimization problem (\ref{eq:tikh_classic}). Because $\cR(A^\ast)=\cR\left(\basa^{\frac{1}{2}}\right)$, this means that all approximate solutions fulfil a source condition (\ref{eq:sc}) with $\mu=\frac{1}{2}$.

Next, note that $\overline{\cR(A^\ast)}=\cN(A)^\perp$. In view of the preceding argument, this yields $\xad,\xa \in\cN(A)^\perp$. Because $\cN(A)^\perp\cap \cN(A)=\{0\}$, Tikhonov regularization automatically takes care of non-injectivity of the forward operator by disallowing any non-zero components in the null-space $\cN(A)$. This also explains the well-known fact that the Tikhonov-approximations converge to the minimum norm solution to (\ref{eq:problem}). 

Later we will often make use of a more specific property of $\xad$ that follows from Proposition \ref{thm:bas}. The source element of $\xad$ is easily available even in practical computations, since we can write
\begin{equation}\label{eq:sc_xad}
\xad=A^\ast \xiad, \quad \xiad=\frac{A\xad-y^\delta}{\alpha}, \quad \mbox{and}\quad \xa=A^\ast \xia, \quad \xia=\frac{A\xa-y}{\alpha}
\end{equation}

For later reference we recall the following result from \cite{GerRamres}.
\begin{proposition}\label{thm:alpharateres}
Let $0<\mu+\nu<1$, then the Tikhonov-approximations for noise-free data $\xa$ satisfy
\begin{equation}\label{eq:alpharate}
\|\basa^{\nu+\mu}\xa-\basa^{\nu+\mu}\xd\|=\bigo (\alpha^{\mu+\nu})
\end{equation}
if and only if (\ref{eq:sc_sum}) holds.
\end{proposition}

To summarize, we suggest to move away from the the standard interpretation that Tikhonov regularization enforces a small norm $\|\xad\|$ of the regularized solutions through minimizing the functional (\ref{eq:tikh_classic}). 
While this is certainly true, we propose the following, stronger, interpretation: Tikhonov regularization is the method of approximating the true solution $\xd$ in $\cR(\basa^\frac{1}{2})$ by controlling the norm of the source element $\xi_\alpha^{\cdot}=\frac{A\xad-y^\cdot}{\alpha}$ in (\ref{eq:sc_xad}) through the regularization parameter $\alpha$. In the remainder of the paper, we collect consequences of this point of view. While many of these are known, the new interpretation may give a better understanding of how regularization works.

We finish this section with an example that shows that the impact of Proposition \ref{thm:bas} is not just theoretical, but immediately affects numerical computations. We consider a Fredholm integral equation of the first kind, 
\begin{equation}\label{eq:exrange}
y(s)=[Ax](s)=\int_0^1k(s,t)x(t)\,dt,\qquad 0\leq s\leq 1,
\end{equation}
with the kernel function $k(s,t)=\left\{
\begin{array}{ll}
s(t-1), & s<t\\ t(s-1), & s\geq t \\
\end{array}
\right.$. With $x^\dag(t)=t$ this is the default setting of the \texttt{deriv2}-example in the RegularizationTools toolbox \cite{regtools}, which we use for the numerical experiments.

Let $x^{\ast\ast}(t)$ be a second primitive of $x$, i.e., $\frac{ d^2 x^{\ast\ast}(t)}{dt^2}=x(t)$. Then one can show by twice partial integration of (\ref{eq:exrange}) that
\[
y(s)=(s-1)x^{\ast\ast}(0)-sx^{\ast\ast}(1)+x^{\ast\ast}(s).
\]
In particular, one easily sees $y(0)=y(1)=0$, hence 
\[
\cR(A)\subseteq \{y\in H^2[0,1]: y(0)=y(1)=0\}.
\]
Since in this particular example $A:L^2[0,1]\rightarrow L^2[0,1]$ is self-adjoint we have $\cR(A)=\cR(A^\ast)$. Because $\xad\in\cR(A^\ast)$ according to (\ref{eq:tikh_sourcerep1}), this means that the Tikhonov-approximated solutions fulfil $\xad(0)=\xad(1)=0$ and are twice differentiable. Again, we mention this is much stronger than the condition $\xad\in L^2[0,1]$ implied by the optimization problem (\ref{eq:tikh_classic}).

\begin{figure}
\includegraphics[width=\linewidth]{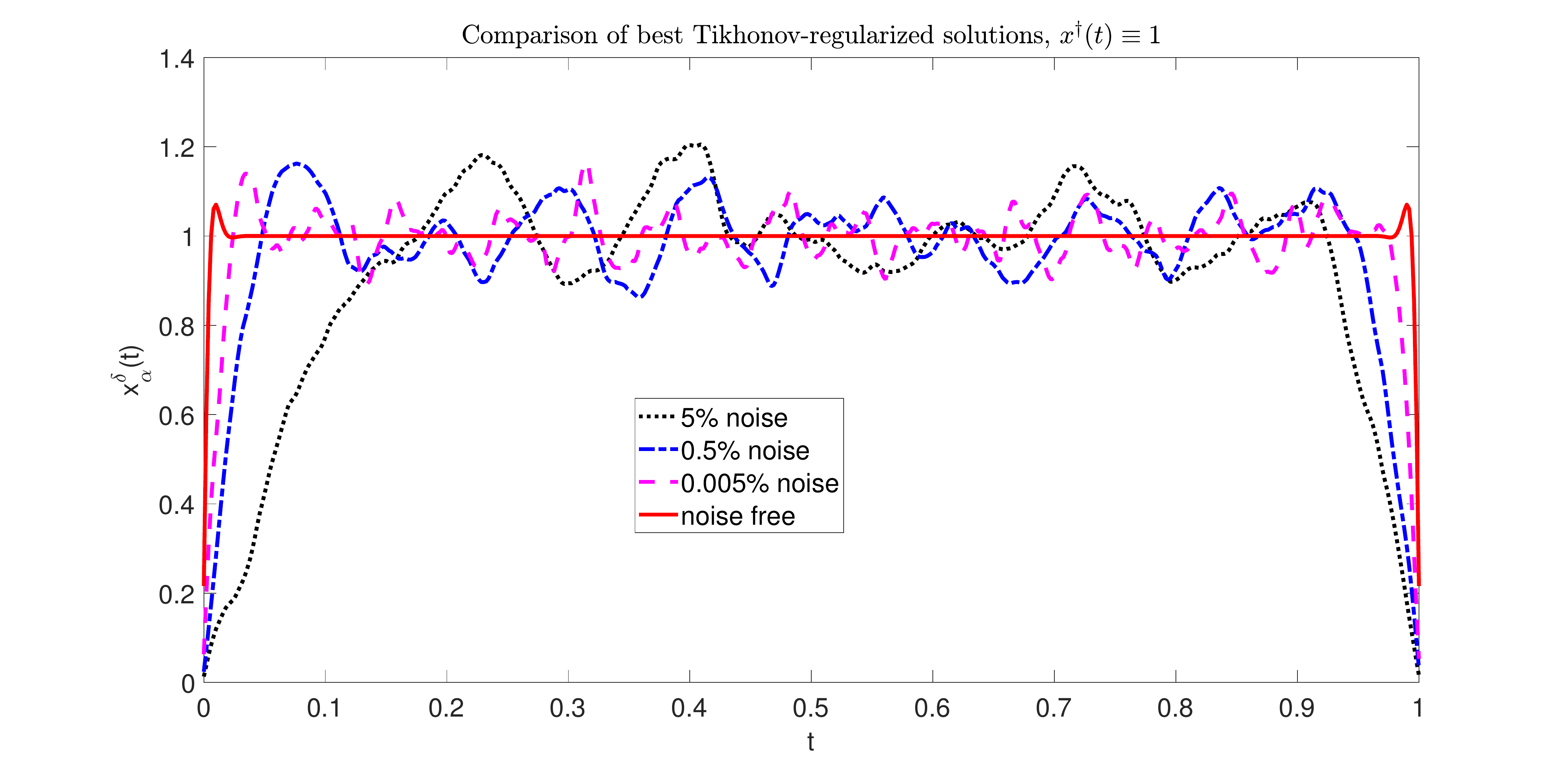}\\\includegraphics[width=\linewidth]{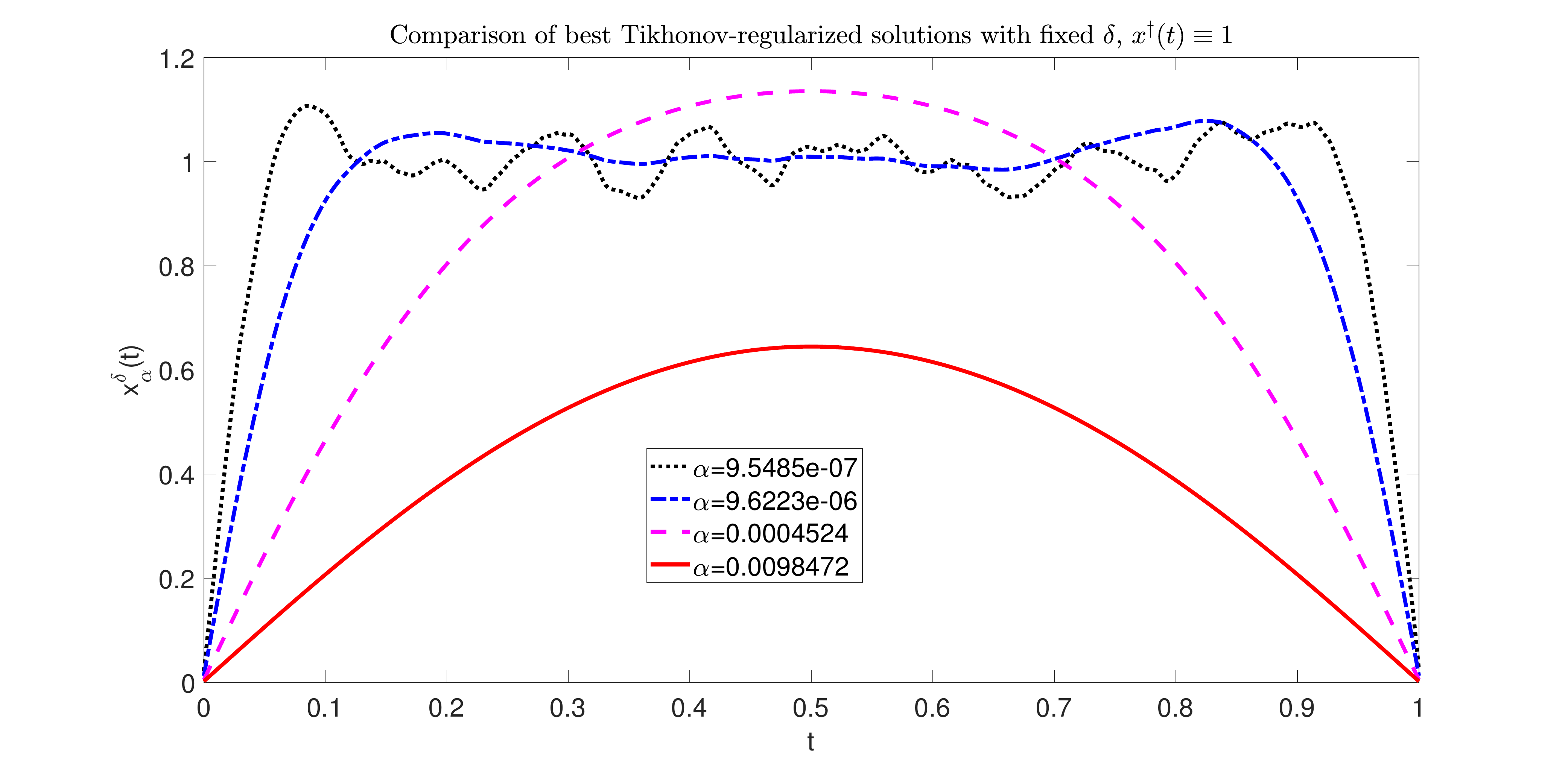}\caption{Tikhonov regularized solutions for (\ref{eq:exrange}) with $\xd\equiv 1$. Top: optimal regularization parameters for various noise levels. Bottom: several regularization parameters for fixed noise level. In all cases, the conditions $\xad(0)=\xad(1)= 0$ implied by $\xad\in\cR(A^\ast)$ can be easily spotted in the reconstructions.}\label{fig:sourceexample}
\end{figure}

The practical effect is that all approximated solutions $\xad$ will be at least close to zero at the boundary points. We can see these implicitly forced boundary conditions clearly in the numerical examples which we sketch in the following and present in an example in Figure \ref{fig:sourceexample}. We use the exact implementation of the forward operator $A$ from the Regularization Tools \cite{regtools}, but for better visualization we set $x^\dag(t) \equiv 1$, calculating $y=Ax$. Now we conduct two experiments. First we add noise to the data, $y^\delta=y+\epsilon$ where $\epsilon$ is Gaussian with zero mean and unit variance (\texttt{randn} in MATLAB). Then we rescale the noise such that $\frac{\|y-y^\delta\|}{\|y\|}=\delta$, where $\delta\cdot100$ is the percentage of the relative error in the data. For $\delta\in\{0,0.0005,0.005,0.05\}$ we calculate the optimal regularization parameter by minimizing $\|\xad-\xd\|$ over 300 values of $\alpha$ between $10^{-10}$ and 0. As result we obtain, as expected, that the solutions $\xad$ approach $\xd$ as $\delta$ decreases. More important, and the purpose of this experiment, is another observation: for all $\delta$ (even noise free), we see that $\xad(0)\approx\xad(1)\approx 0$. This is a property directly inherited from the range restrictions of $A$ and $A^\ast$, as noted above. We see the same effect in the second example, where for fixed $\delta=0.005$ we show four Tikhonov-regularized solutions for different values of $\alpha$. We see how lower values of $\alpha$ at first approximate $\xad\equiv 1$ better in the middle of the interval $(0,1)$ only to start amplifying the noise eventually. However, independent of $\alpha$, $\xad(0)\approx\xad(1)\approx 0$, which illustrates our point.

\section{New convergence theory based on approximate source conditions}\label{sec:appsc}
The core idea of our new approach is to approximate several objects from a common source space. In regularization theory, there is already a smoothness condition that is closely tied to this principle, namely \textit{approximate source conditions} (ASC). This concept has, after gaining some attention in the 2000's, mostly fallen out of favour with the rising popularity of \textit{variational source conditions} (see, e.g., \cite{jens12,jens13,HKPS,HofMat12,hw}), which will not be considered further in this work. Note that many smoothness conditions including source conditions (\ref{eq:sc}) and (\ref{eq:sc_sum}) can be related, and often they are equivalent, see, e.g., \cite{jens12,jens13,GerthKindermann}. However, for our purpose ASCs are ideal smoothness descriptions, and we will discuss and significantly extend existing results in this section, ultimately leading to a new approach to derive convergence rates and parameter choice rules.

The idea behind ASCs is to measure how well an element $\xd\in X$ can be approximated in $\cR(\varphi(\asa))$ with some index function $\varphi: \R\rightarrow \R$, called the benchmark smoothness. More precisely, we define 
\begin{equation}\label{eq:approx_sc}
d_\varphi(R)=\inf_{\|\xi\|\leq R} \|\xd-\varphi\basa\xi\|.
\end{equation}
This concept was developed in \cite{Baumeister} and discussed for Tikhonov regularization in, e.g., \cite{Dana2,Bernd06,Dana}. It was later generalized to other regularization methods, to Banach space regularization and to nonlinear operators, see, e.g., \cite{jens12,jens13,Hein,HeinHof}.

Our renewed interest in approximate source conditions is due to their connection to Proposition \ref{thm:bas}. Under this new paradigm, we know that the \textit{regularized solutions} fulfil a fixed source condition $\xad\in\cR(A^\ast)$ unconditionally. Therefore, Tikhonov regularization can be described as the task of approximating $\xd$ in $\cR(A^\ast)$. This is precisely the concept of the approximate source condition (\ref{eq:approx_sc}) with $\varphi\basa=\basa^\frac{1}{2}$. Therefore, approximate source conditions with benchmark function $\varphi(t)=t^{\frac{1}{2}}$ are a natural way of expressing solution smoothness for Tikhonov regularization. This observation removes a crucial degree of freedom in the formulation of the approximate source condition, namely the function $\varphi$. Instead of being chosen and thus rather arbitrary, it is now fixed by the regularization method. Note that in the formulation of the approximate source condition (\ref{eq:approx_sc}), the source element $\xi$ is not specified further, while  for Tikhonov regularization it is restricted to the form (\ref{eq:sc_xad}).

\subsection{Approximate source conditions under operator powers}
One reason why approximate source conditions had fallen out of favour is because they are, from a mostly theoretical point of view, equivalent to source conditions. This is not surprising, since it appears natural that the approximation rate (\ref{eq:approx_sc}) must be strongly connected to the smoothness of $\xd$. We will discuss further below why ASCs are more powerful than source conditions. For now we note that in \cite[Corollary 3.3]{Dana2} a one-to-one correspondence between a distance function $d_\varphi(R)\leq CR^{-\frac{\mu}{\kappa-\mu}}$, with $\varphi(t)=t^\kappa$, $0<\mu<\kappa$, and a source condition $\xd\in \cR(\basa^\mu)$ was shown. We generalize this in the following, showing that  the source condition is not only equivalent to the approximation rate of $\xd$ in $\cR(\varphi(\asa))$ (\ref{eq:approx_sc}), but also to the approximation rate of operator powers $\basa^\nu \xd$, for certain $\nu\geq 0$ and benchmark smoothness $\varphi(t)=t^\kappa$. 
\begin{theorem}\label{thm:conv_approx}
It is
\begin{equation}\label{eq:gen_approx}
d^\nu_\kappa(R)=\inf_{\|\xi\|\leq R} \|\basa^\nu\xd-\basa^{\nu+\kappa}\xi\|=\bigo( R^\frac{\mu+\nu}{\mu-\kappa})
\end{equation}
for $-\nu<\mu<\kappa$ if and only if
\[
\sum_{n=k}^\infty \langle \xd,v_n\rangle^2=\bigo(\sigma_k^{4\mu}),
\]
implying (\ref{eq:sc_open}). 
\end{theorem}
\begin{proof}
The proof follows the lines of \cite[Theorem 1]{Dana}.
The Lagrange functional for the optimization problem in (\ref{eq:gen_approx}) reads
\[
\|\basa^\nu\xd-\basa^{\nu+\kappa}\xi\|^2+\lambda(\|\xi\|^2-R^2)
\]
with the Lagrange parameter $\lambda$. From the theory of Lagrange multipliers we obtain $\|\xi\|^2=R^2$ and $\xi=(\basa^{2(\nu+\kappa)}+\lambda I)^{-1}\basa^{2\nu+\kappa}\xd$, so that
\begin{equation}\label{eq:d_prep}
d_\kappa^\nu(R)=\lambda(R)\|(\basa)^{2(\nu+\kappa)}+\lambda(R) I)^{-1} \basa^\nu \xd\|.
\end{equation}
Hence, we have to solve
\[
R^2=\|(\basa^{2(\nu+\kappa)}+\lambda I)^{-1}\basa^{2\nu+\kappa}\xd\|^2
\]
for $\lambda$. Using the singular system, this reads
\begin{equation*}
R^2=\sum_{i=1}^\infty \frac{\sigma_i^{8\nu+4\kappa}}{(\sigma_i^{4\nu+4\kappa}+\lambda)^2}\langle \xd,v_i\rangle^2=\lambda^{-2}\sum_{i=1}^\infty \frac{\sigma_i^{8\nu+4\kappa}\lambda^2}{(\sigma_i^{4\nu+4\kappa}+\lambda)^2}\langle \xd,v_i\rangle^2.
\end{equation*}
We apply \cite[Lemma 1]{GerRamres} with $q=8\nu+4\kappa$ and $p=4\nu+4\kappa$ (which is applicable as long as $\mu<\kappa$) and obtain
\[
R^2=\bigo(\lambda^{-2} \lambda^{\frac{8\nu+4\kappa+4\mu}{4\nu+4\kappa}})=\bigo(\lambda^{\frac{4\mu-4\kappa}{4\nu+4\kappa}}),
\]
i.e., $\lambda(R)=\bigo(R^{\frac{8\nu+8\kappa}{4\mu-4\kappa}})$.
With this $\lambda=\lambda(R)$ we then have to evaluate (\ref{eq:d_prep}), which we can write as
\begin{equation*}
d_\kappa^\nu(R)^2=\sum_{i=1}^\infty \frac{\sigma_i^{4\nu}\lambda^2}{(\sigma_i^{4\nu+4\kappa}+\lambda(R))^2}\langle \xd,v_i\rangle^2=\bigo(\lambda(R)^{\frac{4\nu+4\mu}{4\nu+4\kappa}}),
\end{equation*}
where we used \cite[Lemma 1]{GerRamres} with $q=4\nu$ and $p=4\nu+4\kappa$.
Inserting $\lambda(R)$ and taking the square root yields the claim.
\end{proof}

\begin{figure}
\begin{tabular}{c c}
\includegraphics[width=0.45\linewidth]{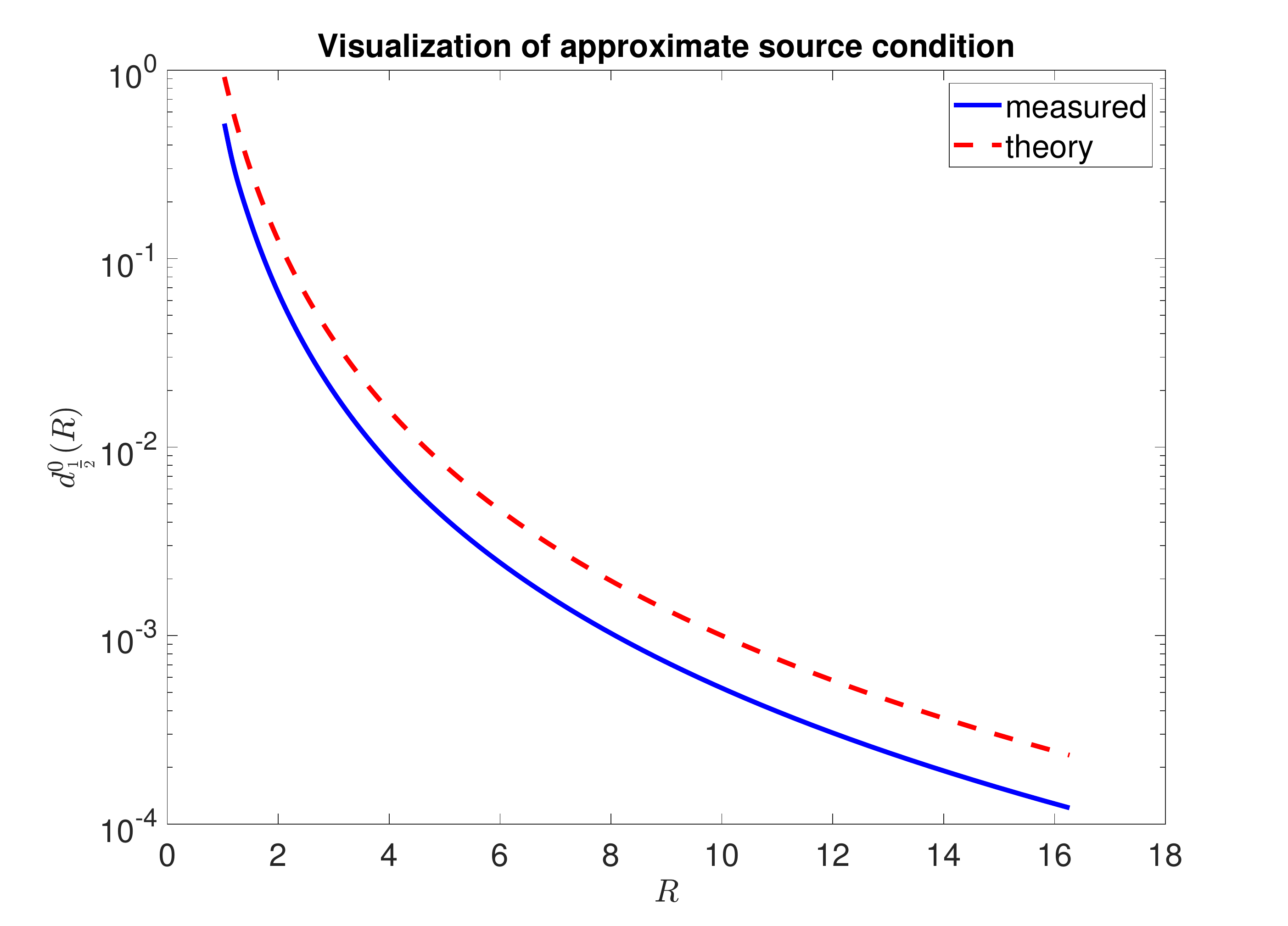}&\includegraphics[width=0.45\linewidth]{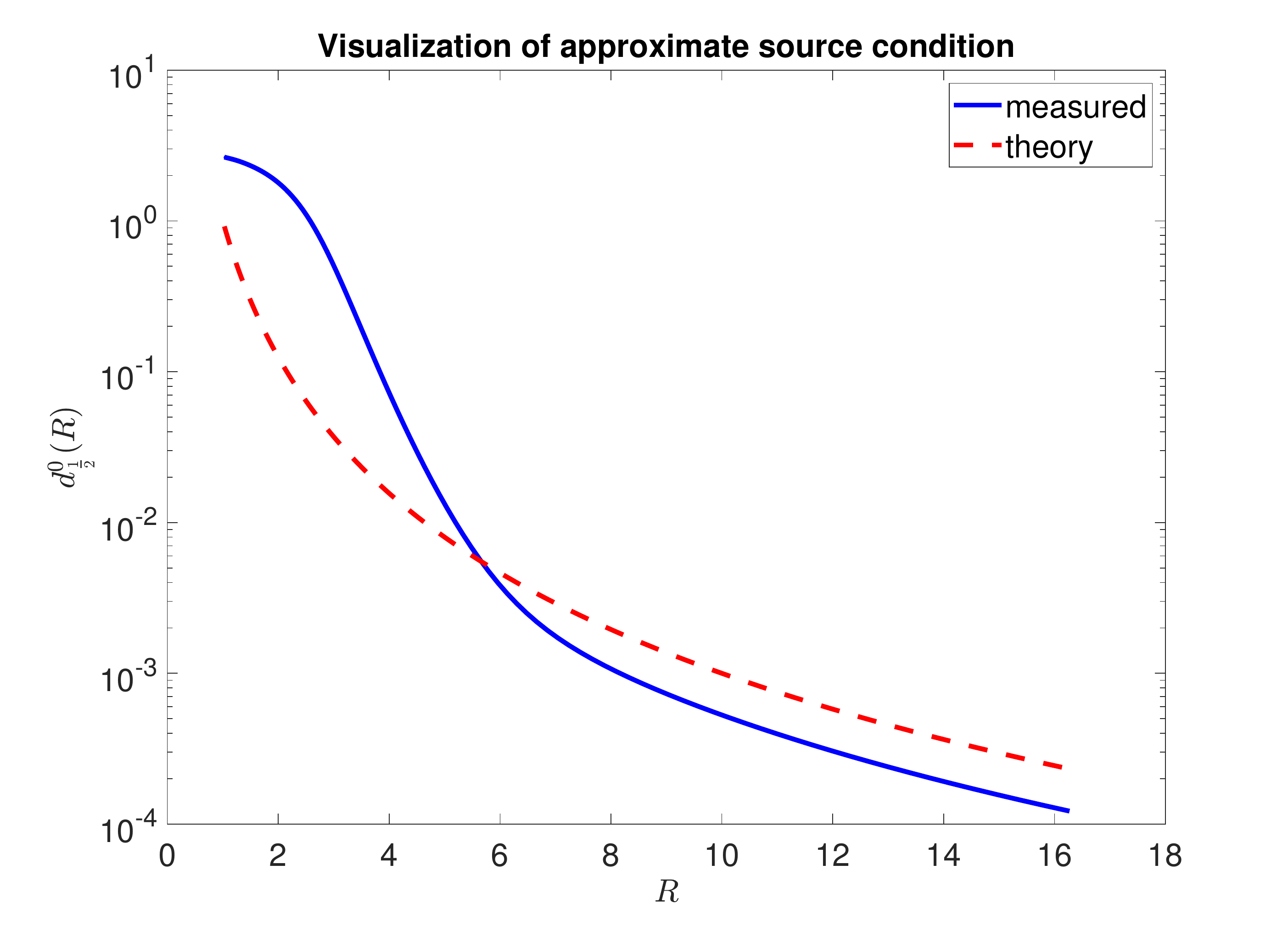}\\
\includegraphics[width=0.45\linewidth]{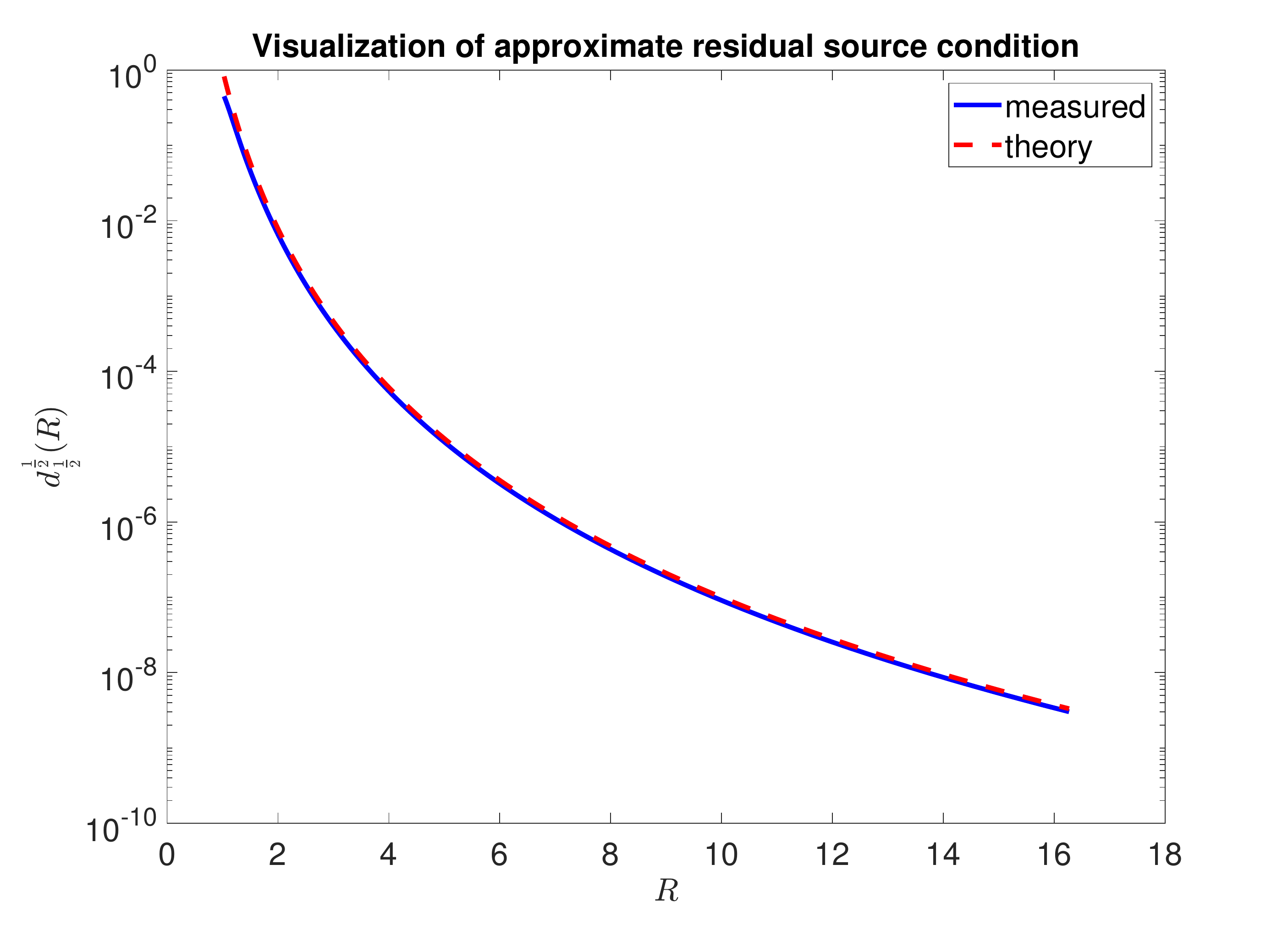}&\includegraphics[width=0.45\linewidth]{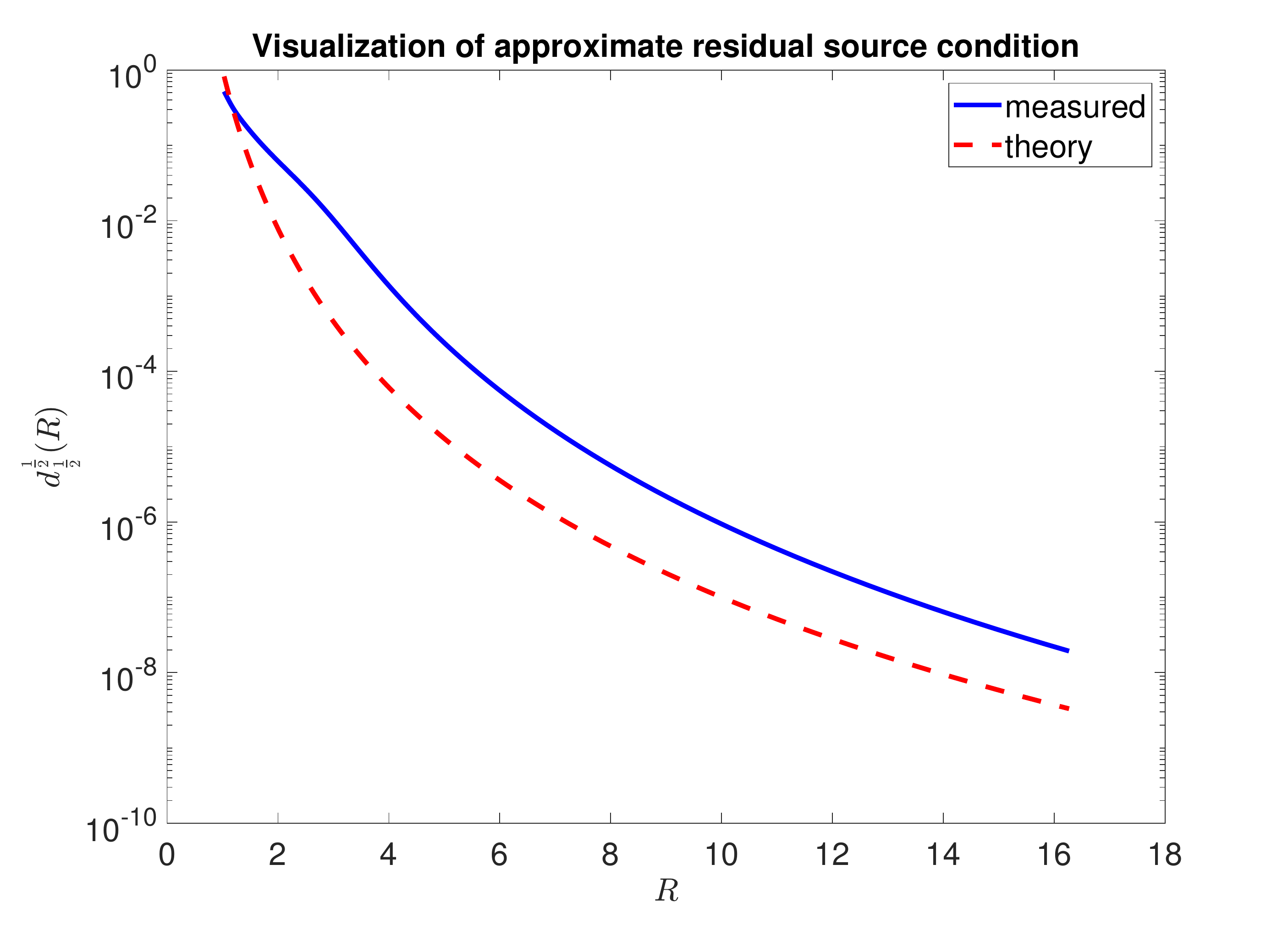}
\end{tabular}\caption{Approximate source condition with $\varphi(t)=t^{\frac{1}{2}}$ for two solutions satisfying (\ref{eq:sc_sum}) with $\mu=0.375$. Left column: (\ref{eq:sc_sum}) holds for all $n>0$. Right column: first 8 components $\langle x^\dag,v_i\rangle$ set to one, (\ref{eq:sc_sum}) holds for all $n>8$. Top row: approximate source condition (\ref{eq:approx_sc}), bottom row: approximate residual source condition, (\ref{eq:gen_approx}) with $\nu=\frac{1}{2}$. Theoretical values in red, measured ones in blue. Both solutions have, up to constants, the same asymptotic for large $R$, but differ notably for small $R$.}\label{fig:appsc}
\end{figure}

The theorem states that eventually, for sufficiently large $R$, approximate source conditions (even under operators, although with some loss of smoothness information) and source conditions are equivalent. ASCs, however, are more precise than source conditions in that they are often able to clearly distinguish between solutions for which the same source condition (\ref{eq:sc_open}) holds, as $d_\varphi(R)$ behaves differently for lower values of $R$, when the asymptotics are not active yet. To see this, consider $\xd=\sum_{i=1}^\infty \langle \xd,v_i\rangle v_i$ fulfilling (\ref{eq:sc_open}) with some $\mu>0$. Any modification of finitely many coefficients  $\langle x^\dag,v_i\rangle$ will change the function (\ref{eq:approx_sc}) while the asymptotic behaviour stays the same. We illustrate this in a numerical example where $A$ is a diagonal operator and the coefficients $\langle \xd,v_i\rangle$ are chosen such that $\xd$ fulfils (\ref{eq:sc_sum}) for chosen $\mu>0$, \cite[Model Problem]{GerRamres}. In the first case, (\ref{eq:sc_sum}) holds for all $n>0$, whereas in the second one (\ref{eq:sc_sum}) holds for all $n>8$, as we set $\langle \xd,v_i\rangle=1$ for $i=1,\dots,8$. The results are shown in Figure \ref{fig:appsc}.

Approximate source conditions have been used to derive convergence rates for Tikhonov regularization before, and we quickly review the ideas. For simplicity, we only recall here the basic result from \cite{Bernd06}. Using the bias-error decomposition 
\begin{equation}\label{eq:decomp}
\|\xad-\xd\|\leq \|\xa-\xd\|+\|\xad-\xa\|
\end{equation}
it was shown that the approximate source condition (\ref{eq:approx_sc}) with $\varphi(t)=t^{\frac{1}{2}}$, which then takes the form $d_\varphi(R)=d_0^{\frac{1}{2}}(R)=\bigo(R^{\frac{\mu}{\mu-\frac{1}{2}}})$ with $d_0^{\frac{1}{2}}(R)$ from Theorem \ref{thm:conv_approx}, implies a bias
\begin{equation*}\label{eq:approxsc_bias}
\|\xa-\xd\|\leq d_\varphi(R)+\sqrt{\alpha}R
\end{equation*}
for Tikhonov regularization. The authors then chose $R=R(\alpha)=\bigo(\alpha^{\mu-\frac{1}{2}})$, which yields 
\[
\|\xa-\xd\|\leq C\alpha^\mu
\]
with a constant $C>0$ that we do not specify further here. Inserting this together with the well-known estimate $\|\xad-\xa\|\leq \frac{\delta}{\sqrt{\alpha}}$ for the noise amplification, (\ref{eq:decomp}) becomes
\[
\|\xad-\xd\|\leq C\alpha^\mu +\frac{\delta}{\sqrt{\alpha}}.
\]
Minimizing the  right-hand size over $\alpha>0$ yields the optimal parameter choice $\alpha=\bigo(\delta^{\frac{1}{2\mu+1}})$ and the optimal convergence rate
\[
\|\xad-\xd\|\leq C \delta^{\frac{2\mu}{2\mu+1}}.
\]
Due to our novel approach, we can show that $R=\bigo(R^{\frac{\mu}{\mu-\frac{1}{2}}})$ is not so much a choice, but a direct consequence of properties of Tikhonov regularization. Assume a source condition (\ref{eq:sc_sum}) with $0<\mu<\frac{1}{2}$. Due to (\ref{eq:alpharate}), it is \[\|A\xa-A\xd\|=\bigo(\alpha^{\mu+\frac{1}{2}}).\] Hence, the source element in $\xa=A^\ast\xia$, $\xia=\frac{A\xa-A\xd}{\alpha}$ (see (\ref{eq:tikh_sourcerep1})) has norm $\|\xia\|=\bigo(\alpha^{\mu-\frac{1}{2}})$. This means Tikhonov regularization approximates $\xa$ in $\cR(A^\ast)$ with a source element of magnitude $\bigo(\alpha^{\mu-\frac{1}{2}})$. By construction, the approximate source condition for $d_0^{\frac{1}{2}}$ (\ref{eq:gen_approx}) with $R=\bigo(\alpha^{\mu-\frac{1}{2}})$ now yields the infimal distance between $\xa$ and $\xd$,
\[
d_0^{\frac{1}{2}}(R)=\inf_{\|\xi\|\leq R} \|\xd-A^\ast\xi\|=\bigo\left( R^\frac{\mu}{\mu-\frac{1}{2}}\right)=\bigo(\alpha^\mu).
\]
Hence we find the well-known bias estimate
\[
\|\xa-\xd\|= \bigo(\alpha^\mu).
\]
From here one could proceed as before and deduct convergence rates via (\ref{eq:decomp}). In the following we demonstrate a new approach that does not require the bias-error decomposition. Instead, it is based on combining two approximate source conditions from Theorem \ref{thm:conv_approx}, one for $\|A\xa-A\xd\|$ and one for $\|\xa-\xd\|$. 

\subsection{A convergence rate prototype based on two ASCs}

Consider two functions $d_\kappa^\nu(R)$ from (\ref{eq:gen_approx}) for fixed $\kappa=\frac{1}{2}$, one with $\nu=0$ and the other with $\nu=\frac{1}{2}$. Let $R$ be fixed. Then the distance functions describe the minimal approximation error ($d_{\frac{1}{2}}^0(R)$) and the minimal residual ($d_{\frac{1}{2}}^{\frac{1}{2}}(R)$) when approximating $\xd$ in $\cR(A^\ast)$, $\xd=A^\ast \xi$, $\|\xi\|\leq R$. With slight abuse of notation (replacing $\bigo$ in (\ref{eq:gen_approx}) by a generic constant $C>0$), we write $d_{\frac{1}{2}}^{\frac{1}{2}}(R)=C R^{\frac{\mu+\frac{1}{2}}{\mu-\frac{1}{2}}}$ as 
\[
R=\left(\frac{1}{C}d_{\frac{1}{2}}^{\frac{1}{2}}(R)\right)^{\frac{\mu-\frac{1}{2}}{\mu+\frac{1}{2}}},
\]
and substitute this $R$ in the ASC $d_{\frac{1}{2}}^0(R)$, which yields
\[
d_{\frac{1}{2}}^0(R)=d_{\frac{1}{2}}^0\left(\left(\frac{1}{C}d_{\frac{1}{2}}^{\frac{1}{2}}(R)\right)^{\frac{\mu-\frac{1}{2}}{\mu+\frac{1}{2}}} \right).
\]
Noting now that $d_{\frac{1}{2}}^0(R)$ is an approximation error $\|\xd-A^\ast  \xi\|$ and $d_{\frac{1}{2}}^{\frac{1}{2}}(R)$ the residual $\|A\xd-AA^\ast \xi\|$, this means (going back to the $\bigo$-notation)
\begin{equation}\label{eq:aprox_rate1}
\|\xd-A^\ast  \xi\|=\bigo\left(\left(\|A\xd-AA^\ast \xi\|^{\frac{\mu-\frac{1}{2}}{\mu+\frac{1}{2}}}\right)^{\frac{\mu}{\mu-\frac{1}{2}}}   \right)=\bigo\left( \|A\xd-AA^\ast \xi\|^{\frac{\mu}{\mu+\frac{1}{2}}} \right)
\end{equation}
for $0<\mu\leq \frac{1}{2}$. On the right-hand side we now find precisely the function $\varphi(t)=t^{\frac{2\mu}{2\mu+1}}$ that characterizes the well-known worst-chase error estimate
\begin{equation}\label{eq:worsterr}
\sup \{\|x\|: x\mbox{ fulfils (\ref{eq:sc}) }, \|Ax\|\leq \delta\}\leq c \delta^{\frac{2\mu}{2\mu+1}}, 
\end{equation}
that holds for any regularization method, see, e.g., \cite{EHN}. Both this and (\ref{eq:aprox_rate1}) are based on the geometry of the underlying spaces, independent of the regularization method. From the point of an approximation problem, the convergence rate is shaped by the way $\xd$ and $A\xd$ can be approxiamted in the same source element space $\mbox{dom}(A^\ast)$. More general, by repeating the previous steps for two generalized ASCs (\ref{eq:gen_approx}), one with $\nu=0$ and $0<\mu<\kappa$, and the other with some $\nu>0$ and the same $\kappa$, we obtain the following.

\begin{theorem}\label{thm:rate_appprox_gen}
Let $\xi\in Y$ and $\nu\geq 0$. Then
\begin{equation*}\label{eq:aprox_rate_gen}
 \|\xd-\basa^\kappa   \xi\|
 =\bigo\left( \|\basa^\nu\xd-\basa^{\nu+\kappa} \xi\|^{\frac{\mu}{\mu+\nu}} \right).
 \end{equation*}
 if and only if $\xd$ satisfies (\ref{eq:sc_open}) with $0<\mu<\kappa$. The corresponding source element $\xi$ has the norm
\begin{equation}\label{eq:xi}
\|\xi\|=\bigo\left([d_\kappa^\nu]^{-1}(\|\basa^\nu\xd-\basa^{\nu+\kappa} \xi\|) \right).
\end{equation}
\end{theorem}
The theorem can be generalized further to approximate source conditions with general benchmark functions $\varphi$ in (\ref{eq:approx_sc}), when also a corresponding approximate residual source condition is assumed. However, this is not the scope of this paper. Note that earlier, with $\nu=\kappa=\frac{1}{2}$ we were only able to show the optimal convergence order if $\xd$ satisfies (\ref{eq:sc_sum}) for $0<\mu<\frac{1}{2}$. To extend this to $0<\mu<1$, we have to consider the alternative source representation (\ref{eq:tikh_sourcerep2}) which can be interpreted as approximating $\xd$ in $\cR(A^\ast A)$. Now we can apply Theorem \ref{thm:rate_appprox_gen} with $\kappa=1$ and $\nu=\frac{1}{2}$ to obtain (\ref{eq:aprox_rate1}) for all $0<\mu\leq 1$, from which we will find the a-priori convergence rates for Tikhonov regularization further down below in this section. Before this, we will show that approximate source conditions can also be used in the case of a well-posed problem, for discrete problems, and to describe the smoothness of the noise.

\subsection{ASCs for well-posed problems}

We start with the case of a well-posed problem.
\begin{lemma}
Let $A:X\rightarrow Y$ be a bounded linear operator with $\cR(A)=\overline{\cR(A)}$ and $x\in X$. Then, there is $R^\ast$ such that, for \[d^\nu_\kappa(R):=\inf_{\|\xi\|\leq R} \|\basa^\nu\xd-\basa^{\nu+\kappa}\xi\|\], it is
\begin{equation}\label{eq:gen_approx_closed}
d^\nu_\kappa(R)= \|\basa^\nu\xd-\basa^\nu P_{\cR(\basa^\kappa)}(\xd)\|. 
\end{equation}
for all $R\geq R^\ast$, where $P_{\cR(\basa^\kappa)}(\cdot)$ are the orthogonal projectors onto $\cR(\basa^\kappa)$. For $R<<R^\ast$ we have 
\begin{equation}\label{eq:gen_approx_closed_incomplete}
d^\nu_\kappa(R)\leq \frac{c}{\sqrt{R}}
\end{equation}
with some constant $c>0$.
\end{lemma}
\begin{proof}
With $\cR(A)=\overline{\cR(A)}$ we have $\cR(A^*)=\overline{\cR(A^*)}$ and \[\cR(\basa^\kappa)=\overline{\cR(\basa^\kappa)}.\] Hence $P_{\cR(\basa^\kappa)}(\xd)$ is uniquely defined, and 
\[
\inf_{\xi} \|\basa^\nu\xd-\basa^{\nu+\kappa}\xi\|=\|\basa^\nu\xd-\basa^\nu P_{\cR(\basa^\kappa)}(\xd)\|.
\]
We can now write $P_{\cR(\basa^\kappa)}(\xd)=\basa^\kappa w$. Setting $R^\ast=\|w\|$ then yields the first part. In particular, if $\xd\in \cR(\basa^\kappa)$, 
\[
\|\basa^\nu\xd-\basa^\nu P_{\cR(\basa^\kappa)}(\xd)\|=0\]
for $R\geq R^\ast$.

For smaller values of $R$, we follow the lines of the proof of Theorem \ref{thm:conv_approx}. The Lagrange functional for (\ref{eq:gen_approx_closed}) reads
\[
\|\basa^{\nu+\kappa}w-\basa^{\nu+\kappa}\xi\|^2+\lambda(\|\xi\|^2-R^2).
\]  
It follows again that $R^2=\|\xi\|^2$, and we can write
\begin{equation}\label{eq:xi_gen_closed}
\xi=(\basa^{2(\nu+\kappa)}+\lambda I)^{-1}\basa^{2(\nu+\kappa)}w=\sum_{i=1}^\infty \frac{\sigma_i^{2(\nu+\kappa)}}{\sigma_i^{2(\nu+\kappa)}+\lambda}\langle w,v_i\rangle.
\end{equation}
The next step is to relate $R$ and $\lambda$. Since $\frac{\sigma_i^{2(\nu+\kappa)}}{\sigma_i^{2(\nu+\kappa)}+\lambda}\leq 1$, $R=\|\xi\|\leq \|w\|$ independent of $\lambda$. Inspecting the sum in (\ref{eq:xi_gen_closed}), we see that this is a good estimate whenever $\lambda$ is sufficiently small, and that $\lim_{\lambda\rightarrow 0} \|\xi\|=\|w\|$. In particular, for $\lambda\rightarrow 0$, $\xi\rightarrow w$ and thus (\ref{eq:gen_approx_closed}) follows.

However, we also have lower bound
\[
\|\xi\|^2=\sum_{i=1}^\infty \frac{\sigma_i^{4(\nu+\kappa)}}{(\sigma_i^{2(\nu+\kappa)}+\lambda)^2}\langle w,v_i\rangle^2\geq \frac{\sigma_0^{4(\nu+\kappa)}}{(\sigma_0^{2(\nu+\kappa)}+\lambda)^2} \|w\|^2
\]
where we used that, since $\cR(A)=\overline{\cR(A)}$, $\sigma_i\geq \sigma_0>0$ for all $i \in \N$. Let now $\lambda >>\sigma_0^{4(\nu+\kappa)}$, then 
\[
\frac{\sigma_0^{4(\nu+\kappa)}}{(\sigma_0^{2(\nu+\kappa)}+\lambda)^2}\approx \frac{\sigma_0^{4(\nu+\kappa)}}{\lambda^2}
\]
and hence $R=\|\xi\|\approx \frac{c\|w\|}{\lambda}$, or $\lambda \approx \frac{c}{R}$. Similar to the proof of Theorem \ref{thm:conv_approx}, we insert the expression for $\xi$ in (\ref{eq:gen_approx_closed}), which yields
\[
d^\nu_\kappa(R)=\lambda\|(\basa^{(2(\nu+\kappa)}+\lambda I)^{-1} \basa^{\nu+\kappa}w\|.
\]
Because $\|(\basa^{(2(\nu+\kappa)}+\lambda I)^{-1} \basa^{\nu+\kappa}\|\leq \frac{1}{\sqrt{\lambda}}$, and, as noted earlier, $\lambda \approx \frac{c}{R}$, we obtain
$d^\nu_\kappa(R)\approx \frac{c}{\sqrt{R}}$, which holds for sufficiently small values of $R$.
\end{proof}
The theorem essentially says that in the well-posed case, there is no difference in approximation rate for $\xd$ or the data $A\xd$. In addition to this, it can be used to understand ill-posedness in the discrete setting. Let $A$ and $\xd$ be the infinite-dimensional operator and solution to (\ref{eq:problem}), and consider their discrete, finite dimensional approximations $A_{mn}\in \R^{m\times n}$ and $\xd_n\in\R^n$ for the discretization levels $m,n\in \N$. $A_{mn}$ is compact and has closed range. One can show that if $A$ is ill-posed, $\xd_n =(A_{mn}^TA_{mn})^\mu \xi_n^\mu$ has a solution for all $\mu>0$, but with $\|\xi_n^\mu\|\rightarrow \infty$ as $n$ and/or $\mu$ go to infinity \cite{RamRei}. In other words, $\xd_n$ fulfils a source condition with respect to $A_{mn}$ for all $\mu\geq 0$, but with exploding source element.

\begin{figure}
\includegraphics[width=\linewidth]{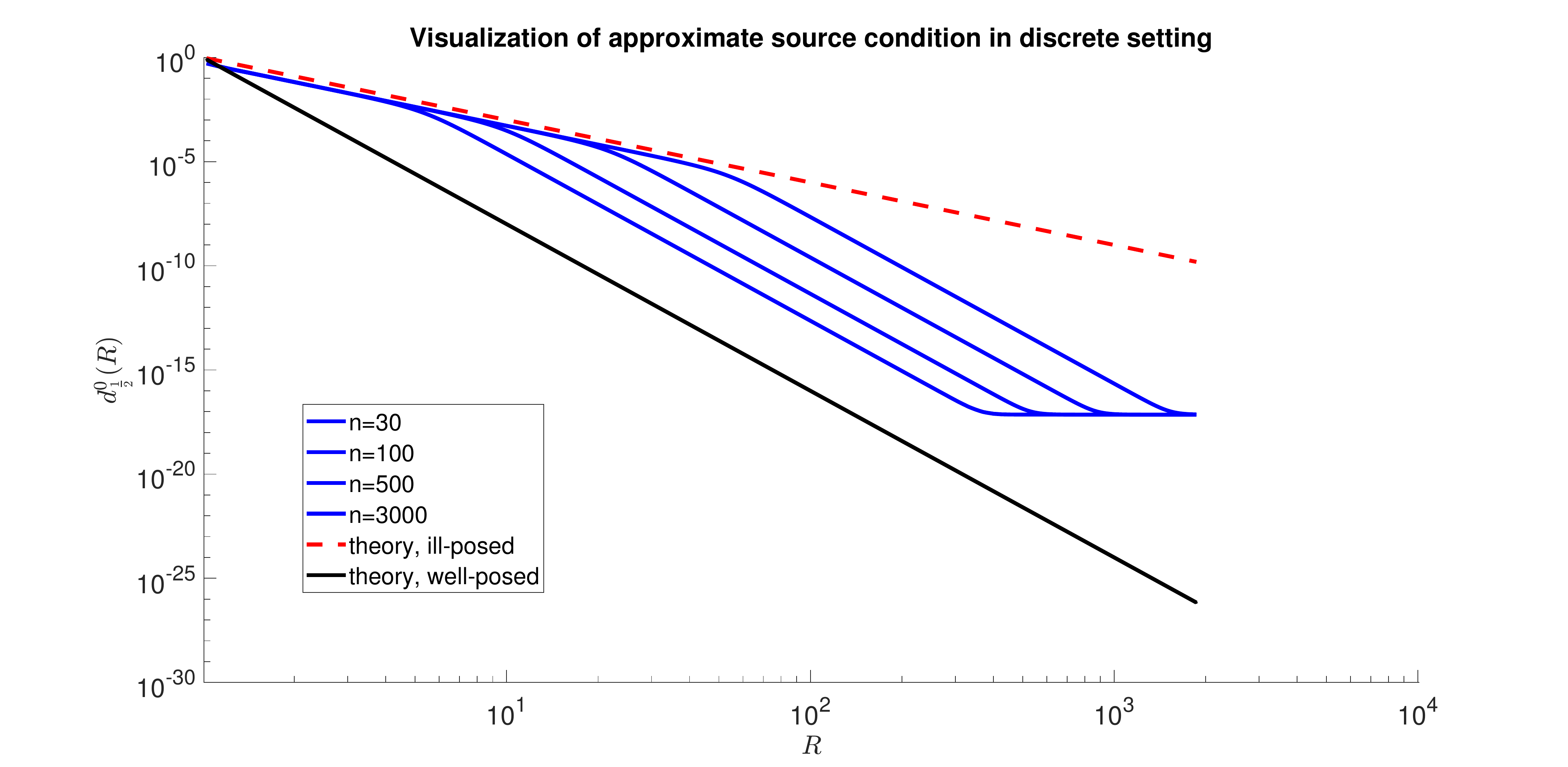}
\caption{Approximate source conditions in the discrete setting. Numerical experiment for source element growth. Diagonal operator and $\xd$ such that a source condition (\ref{eq:sc_open}) holds with $\mu=0.375$, with four different discretization levels. For small $R$, the infinite dimensional approximation properties hold (\ref{eq:gen_approx}). For larger $R$, increasing with the discretization level, the well-posed approximation properties take over. Eventually, as $\xd_n=A^\ast w_n$, $d(R)\approx0$ for $R\geq\|w_n\|$.}\label{fig:ascdiscrete}
\end{figure}

\begin{figure}
\includegraphics[width=\linewidth]{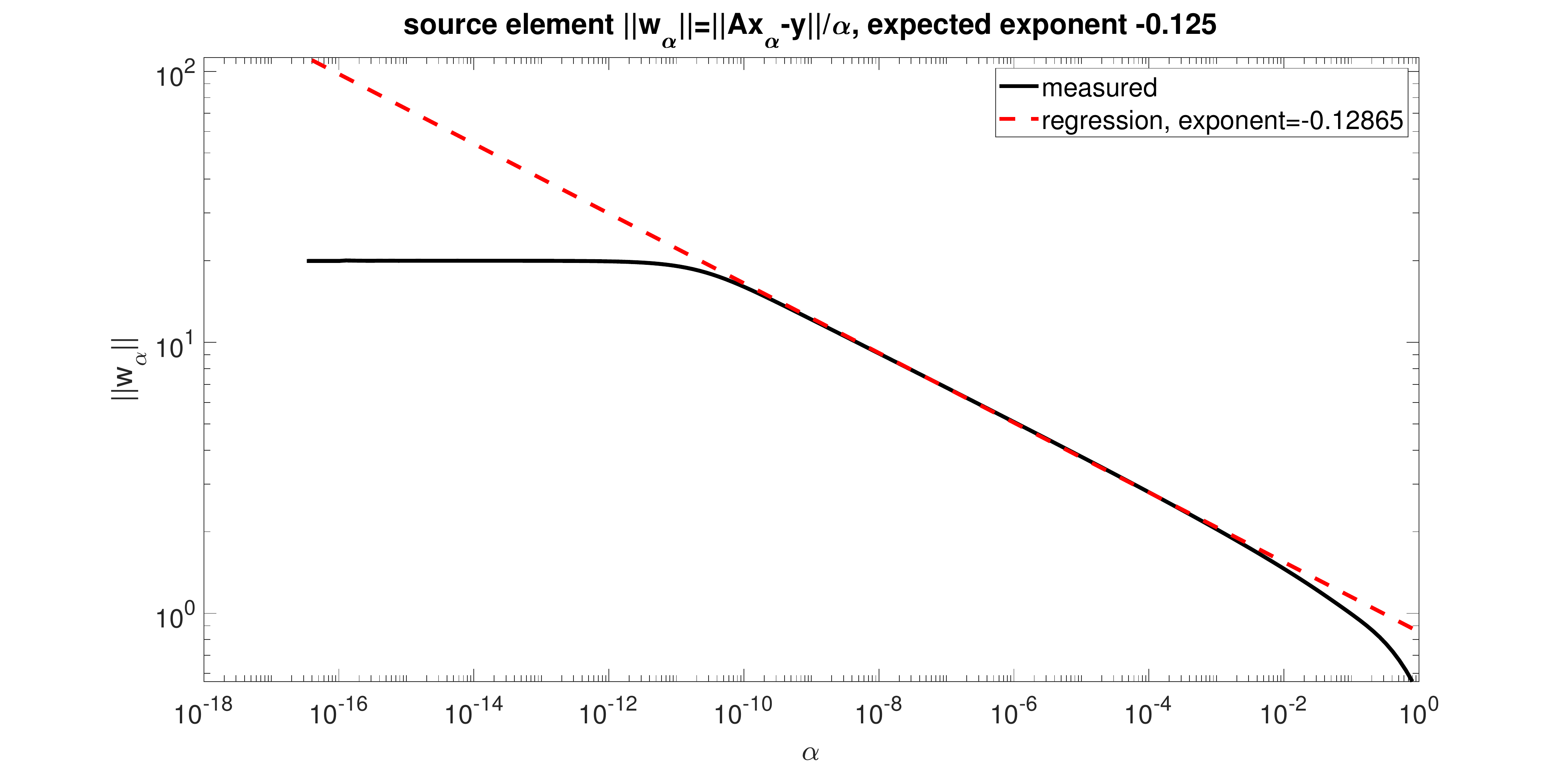}
\caption{Numerical experiment for source element growth. Diagonal operator and $\xd$ such that a source condition (\ref{eq:sc_open}) holds with $\mu=0.375$. Norm of the source element $\xia$ of the Tikhonov-regularized solutions $\xa=A^\ast \xia$ plotted against $\alpha$. For $\alpha$ large enough, we find $\|\xia\|=\bigo(\alpha^{\mu-\frac{1}{2}})=\bigo(\alpha^{-0.125})$ as predicted in the infinite-dimensional setting. Eventually, as $\alpha$ decreases, we have that numerically $\xd_n= A^T_{mn} \xia$, hence $\|\xia\|=\bigo(1)$.}\label{fig:appscwa}
\end{figure}

For Tikhonov regularization, it is most important how $\xd$ can be approximated in $\cR(A^\ast)$ and $\cR(A^\ast A)$. Due to the considerations above, there is $\xi_n\in \R^n$ such that $\xd_n=A_{mn}^\ast \xi_n$, and one would expect to see the approximate source conditions being that of a well-posed situation. This is indeed the case when the norm of the source element is large enough. However, when the source element is small, the approximation properties of $\xd_n$ in $\cR(A^T_{mn})$ are (almost) identical to the ones of the infinite dimensional $\xd$ in $\cR(A^\ast)$. Figure \ref{fig:ascdiscrete} demonstrates this for a numerical example with different discretization levels. For small $R$, the discrete solutions $\xd_n$ follow the theoretical approximate source condition $d_{\frac{1}{2}}^{0}(R)$ (\ref{eq:gen_approx}). After a certain value $R_n$, that increases with increasing discretization level, the solutions follow the ASC in the well-posed sense (see (\ref{eq:gen_approx_closed_incomplete}), $d(R)\approx \frac{c}{\sqrt{R}}$) for $R>R_n$. Eventually, when $R$ is sufficiently large, we have $d(R)\approx 0$, i.e., (\ref{eq:gen_approx_closed}). Due to numerical effects such as round-off errors we do not reach $d(R)=0$. Therefore, approximate source conditions are even useful in a discrete setting, and can be used to determine whether a regularization methods works with the infinite-dimensional approximation rate, or the ``discretization has saturated'', i.e., the approximation follows that of a well-posed problem. We explain this for the example of Tikhonov regularization. The key observation for this is that, due to (\ref{eq:tikh_sourcerep1}), one can always calculate the source element $\xia=\frac{\|A\xa-y\|}{\alpha}$. Due to Proposition \ref{thm:alpharateres}, $\|A\xa-y\|=\bigo(\alpha^{\mu+\frac{1}{2}})$ iff $\xd$ fulfils (\ref{eq:sc_open}) for $0<\mu<\frac{1}{2}$, and $\|A\xa-y\|=\bigo(\alpha)$ if $\xd$ fulfils (\ref{eq:sc_open}) with $\mu\geq \frac{1}{2}$. Hence, $\|\xia\|=\bigo(\alpha^{\mu-\frac{1}{2}})$ in the former case, and $\|\xia\| =\bigo(\alpha/\alpha)=\bigo(1)$ in the latter case which also represents the situation of insufficient discretization. One can observe those two distinct phase by plotting $\|\xia\|$ as function of $\alpha$, see Figure \ref{fig:appscwa}.

\subsection{An ASC for data noise}

Before using Theorem \ref{thm:rate_appprox_gen} to derive convergence rates and parameter choice rules, we need to address the case of noise in the data, which was not under consideration yet in this section. Since we can use ASCs to describe the smoothness of the (unperturbed) residuals, it is natural to investigate if one can do the same with the noise. For the additive noise model (\ref{eq:noise}), we have the following result.

\begin{theorem}
Let $\epsilon=y-y^\delta$, $\|\epsilon\|=\delta$ be the additive noise component in the data to the problem (\ref{eq:problem}). Then, for $\kappa>0$,
\begin{equation}\label{eq:res_gen_approxsc}
d_\epsilon(R)=\inf_{\|\xi\|\leq R} \|\epsilon-A(A^\ast A)^\kappa\xi\|\leq \delta \qquad \mbox{ for all }  R>0.
\end{equation}
\end{theorem}
\begin{proof}
We use again the theory of Lagrange multipliers, i.e., minimize
\[
\|A(A^\ast A)^\kappa\xi -\epsilon\|^2+\lambda(\|\xi\|^2-R^2).
\]
The first order condition is $AA^\ast(AA^\ast\xi-\epsilon)+\lambda\xi=0$, which yields
\[
\xi=((A^\ast A)^{1+2\kappa} +\lambda I)^{-1}(A^\ast A)^{\frac{1}{2}+\kappa} \epsilon
\]
Inserting this into (\ref{eq:res_gen_approxsc}) yields
\[
d(R)=\lambda\|((A^\ast A)^{1+2\kappa}+\lambda I)^{-1}\epsilon\|\leq\lambda\|((A^\ast A)^{1+2\kappa}+\lambda I)^{-1}\|\|\epsilon\|,
\]
and since $\|((A^\ast A)^{1+2\kappa}+\lambda I)^{-1}\|\leq \frac{1}{\lambda}$,
we have 
\[
d(R)\leq \frac{\lambda}{\lambda}\|\epsilon\|=\|\epsilon\|=\delta.
\]
\end{proof}
The point here is that approximate source conditions are non-informative for the noise since they are independent of $R$. While this appears to be rather uninteresting at first sight, it is crucial for the determination of the regularization parameter and thus the convergence rates for Tikhonov regularization.

\subsection{Convergence rates and parameter choice}
\begin{figure}\begin{center}
\includegraphics[width=0.7\linewidth]{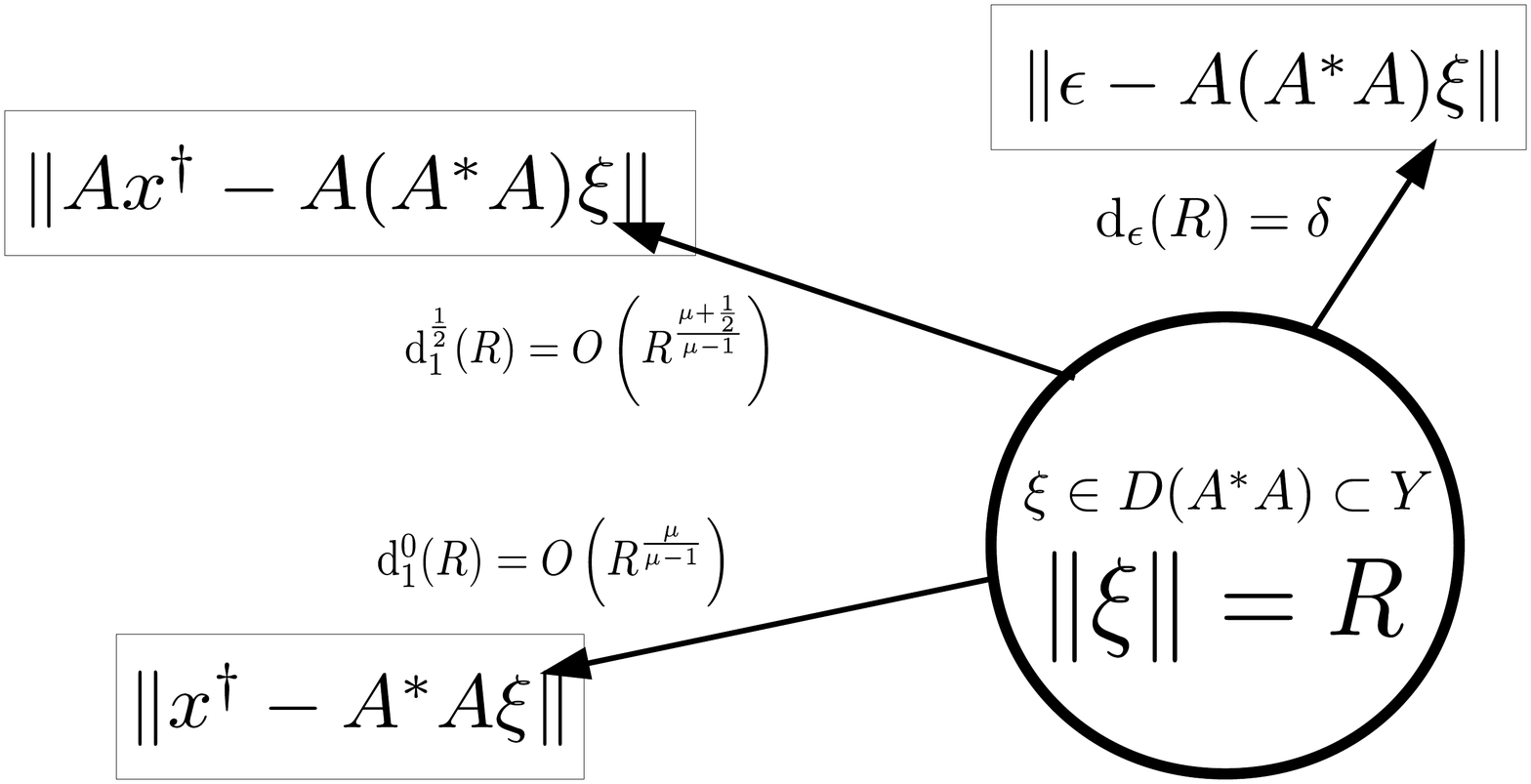}\end{center}\caption{Visualization of the approximation scheme based on $D(A^\ast A)$, i.e., $\kappa=1$. The arrows are labeled with the corresponding approximate source conditions $d_\kappa^\nu$ (\ref{eq:gen_approx}) for $\kappa=1$ and $\nu=\frac{1}{2}$ and $\nu=0$, respectively, and $D_\epsilon$ (\ref{eq:res_gen_approxsc}).}\label{fig:approxoverview}
\end{figure}

Convergence properties of Tikhonov regularization are well studied. Below we discuss these from the apoximation-based point of view. While this does not yield ground breaking new results, it provides valuable insight. The main idea behind the new strategy is that exact solution $\xd$, exact data $y$, and noise $\epsilon$ are all approximated from the same source space, either $D(A^\ast A)$ or $D(A^\ast)$. Figure \ref{fig:approxoverview} visualizes this for the approximation in $D(A^\ast A)$.

It is clear from Theorem \ref{thm:bas} that the role of the regularization parameter is to control the growth of the source element. This can be used to obtain and interpret convergence rates and parameter choices through our novel point of view. Note that in Theorem \ref{thm:rate_appprox_gen} and the discussion preceding it, we have established that the optimal convergence rate is a result of the relation between smoothness of the data in the image space $Y$ and the smoothness of $\xd$ in the pre-image space $X$. Linking two approximate source conditions, one for each object, lead to the relation (\ref{eq:aprox_rate1}), which is the basis for our approach. We may also interpret this as choosing the source element $\xiad$ in $\xad=A^\ast \xiad$ in relation to the residual. What is missing is the link to the regularization method, here Tikhonov regularization. To do this we need to distinguish between the two source representations (\ref{eq:tikh_sourcerep1}) and (\ref{eq:tikh_sourcerep2}) for approximating in $\cR(A^\ast)$ ($\kappa=\frac{1}{2}$ below)  or $\cR(A^\ast A)$ ($\kappa=1$ below), respectively. Depending on $\kappa$, we define the source representation for our approximate $x_\cdot=\basa^\kappa \xi_\kappa$

The first observation is that, due to (\ref{eq:res_gen_approxsc}) and the additive noise model (\ref{eq:noise}), the residuals $\|A\xad-y^\delta\|$ will stagnate around the noise level $\delta$. This effect has also been explained in \cite{GerRamres}. Therefore, in case of noisy data, one should not aim for a residual below the noise level. We apply (\ref{eq:xi}) with $\nu=\frac{1}{2}$ and replace the argument with $\delta$. This yields a source element of magnitude 
\begin{equation}\label{eq:xi_kappa_delta}
\|\xi_\kappa\|=\bigo\left([d_\kappa^{\frac{1}{2}}]^{-1}(\delta) \right)=\bigo\left(\delta^{\frac{\mu-\kappa}{\mu+\frac{1}{2}}} \right).
\end{equation}
In the noise free case, we have from (\ref{eq:tikh_sourcerep1}) and (\ref{eq:tikh_sourcerep2}), respectively, 
\[
\|\xi_\kappa\|=\frac{\basa^{1-\kappa}(\xa-\xd)}{\alpha}=\bigo(\alpha^{\mu-\kappa})
\]
iff $\xd$ satisfies (\ref{eq:sc_open}) with $0<\mu<\kappa\leq 1$. Equating both expressions for $\|\xi_\kappa\|$ yields, ignoring the $\bigo$, $\delta^{\frac{\mu-\kappa}{\mu+\frac{1}{2}}}=\alpha^{\mu-\kappa}$. Solving this for $\alpha$ gives
\begin{equation}\label{eq:a_priori_asc}
\alpha=\bigo\left(\delta^{\frac{1}{\mu+\frac{1}{2}}} \right)=\bigo\left(\delta^{\frac{2}{2\mu+1}} \right).
\end{equation}
This is the classical a-priori result which holds for $0<\mu<1$.  The optimal convergence rate now follows immediately from Theorem \ref{thm:rate_appprox_gen}, since we forced 
\[
 \|\basa^\nu\xd-\basa^{\nu+\kappa} \xi_k\|=\bigo(\delta).
 \]
Alternatively, we may use Theorem \ref{thm:conv_approx} to deduce the convergence rate. Namely, because the a-priori parameter choice yields a source norm $\|\xi_\kappa\|=\bigo\left(\delta^{\frac{\mu-\kappa}{\mu+\frac{1}{2}}} \right)$, we have
\[
\|\xad-\xd\|=\bigo\left(d_{\kappa}^0(\|\xi_\kappa\|) \right)=\bigo\left(\delta^\frac{\mu}{\mu+\frac{1}{2}}\right)=\bigo\left(\delta^\frac{2\mu}{2\mu+1}\right).
\]

Our approach also yields optimal convergence rates for the discrepancy principle, and that the parameter choice obtained via discrepancy principle and the a-priori rate coincide. According to the discrepancy principle, one should choose $\alpha$ such that $\delta<\|A\xad-y^\delta\|\leq \tau \delta$ for some $\tau>1$ (or similarly in closely related formulations such as (\ref{eq:dp1})). Let $\xd$ fulfil (\ref{eq:sc_open}) with $0<\mu<\frac{1}{2}$, and consider the source element for $\xad$ from (\ref{eq:tikh_sourcerep1}). It is $\xad=A^\ast \frac{A\xad-y^\delta}{\alpha}$, and due to the discrepancy principle $\frac{\delta}{\alpha}<\|\xi_{\frac{1}{2}}\|\leq\frac{\tau \delta}{\alpha}$. Equating $\|\xi_\frac{1}{2}\|=\bigo(\frac{\delta}{\alpha})$ with the expression from (\ref{eq:xi_kappa_delta}) with $\kappa=\frac{1}{2}$, we again find $\alpha=\bigo(\delta^{\frac{2}{2\mu+1}})$ and thus $\|\xad-\xd\|=\bigo(\delta^{\frac{2\mu}{2\mu+1}})$. In particular, the discrepancy principle and the a-priori choice can only differ in constants.

\section{Saturation}\label{sec:sat}
A much discussed feature of Tikhonov-regularization is saturation, which means that an arbitrarily high $\mu>$ in a source condition (\ref{eq:sc}) will not lead to the corresponding error bound (\ref{eq:worsterr}), instead the highest convergence rate is $\|\xad-\xd\|=\bigo(\delta^{\frac{2}{3}})$ which holds for all $\mu\geq 1$. It is one of the main disadvantages of Tikhonov-regularization in view of regularization theory. The approximation-based approach based on Theorem \ref{thm:bas} can explain this phenomenon intuitively. Due to (\ref{eq:tikh_sourcerep1}) we know $\xad\in \cR(A^\ast)=\cR(\basa^{\frac{1}{2}})$ with known source element since $\xad=A^\ast\frac{A\xad-y^\delta}{\alpha}$. Now we can distinguish two cases: $\xd$ is less smooth than $\xad$, or $\xd$ is smoother than $\xad$. In the former case, when $\xd$ fulfils a source condition (\ref{eq:sc}) with $\mu<\frac{1}{2}$, it is necessary that $\|\xiad\|:=\frac{\|A\xad-y^\delta\|}{\alpha}\rightarrow \infty$ in order to approximate $\xd$ arbitrarily well. Note that $\xiad$ is driven by the residual. Further, this situation yields $\xad-\xd\in \cR(\basa^\mu)$, and hence we can expect the convergence rate to be of the same order as the worst case error (\ref{eq:worsterr}). This is the unsaturated case. Now let $\xd\in\cR(\basa^\mu)$ with $\mu\geq\frac{1}{2}$. It is well known that the discrepancy principle as parameter choice rule does not yield optimal convergence rates anymore. This is easily seen as the solution smoothness with respect to the residual is fixed, $\xad\in\cR(A^\ast)$, hence $\xad-\xd\in \cR(\basa^{\frac{1}{2}})$, and the convergence rate saturates, $\|\xad-\xd\|\leq c\delta^{\frac{1}{2}}$ ((\ref{eq:worsterr}) with $\mu=\frac{1}{2}$). However, if $\xd \in \cR(\basa^\mu)$ with $\mu>\frac{1}{2}$, then in particular $\xd\in \cR(\asa^{\frac{1}{2}})$, i.e., there is $\xi$ such that $\xd=\basa^{\frac{1}{2}}\xi$. Therefore in this situation the source element of regularized solutions, $\|\xiad\|=\frac{\|A\xad-y^\delta\|}{\alpha}$ must be bounded from below for all $\alpha,\delta$ (as otherwise $\xad=A^\ast\xiad\rightarrow 0$) and above (since otherwise $\xad=A^\ast\xiad$ diverges as $\|\xiad\|\rightarrow \infty$). This means that there must be constants $0<c_1,c_2<\infty$ such that
\begin{equation}\label{eq:whatever}
c_1 \alpha \leq \|A\xad-y^\delta\|\leq c_2\alpha \mbox{ for } 0<\alpha<\alpha_0
\end{equation}
or, in other words, $\|A\xad-y^\delta\|=\bigo(\alpha)$ is necessary for convergence if $\xd$ fulfills a source condition with $\mu>\frac{1}{2}$. This is in accord with Proposition \ref{thm:alpharateres}, which states $\|A\xa-A\xd\|=\bigo(\alpha)$ for this smoothness of $\xd$. Because the residual now bears no more information about the solution smoothness, the discrepancy principle no longer yields the optimal convergence rates.
On the other hand, it is well-known that if $\xd$ fulfils a source condition with $\frac{1}{2}\leq \mu<1$, the a-priori parameter choice $\alpha\sim \delta^{\frac{2}{2\mu+1}}$ still yields convergence rates of optimal order. Also this can be explained through the smoothness of the approximate solutions. The crucial observation is that one may switch to a higher smoothness for the approximate solutions in the following sense. We now consider the representation (\ref{eq:tikh_sourcerep2}). In the absence of noise we see that $\xa\in \cR(\asa)$, and again we simply find the source element in $\xa=\asa \xia$ and its norm, $\|\xia\|=\frac{\|\xa-\xd\|}{\alpha}$. Now the reconstruction error itself is the driver of the solution smoothness, instead of the residual, which leads to an argument similar to the one surrounding (\ref{eq:whatever}). Namely, if $\xd\in \cR(\basa^\mu)$ with $\mu<1$, then $\frac{\|\xad-\xd\|}{\alpha}\rightarrow \infty$ is necessary to approximate $\xd$ in $\cR(\asa)$ arbitrarily well. On the other hand, if $\xd\in \cR(\basa^\mu)$ with $\mu\geq1$, then
\[
\|\xa-\xd\|=\bigo(\alpha)
\]
is necessary for the convergence $\xa\rightarrow \xd$. This is, again, in line with Proposition \ref{thm:alpharateres}, which states $\|\xa-\xd\|=\bigo(\alpha)$ for $\mu\geq 1$. 

It remains to discuss the impact of the noise in the case $\frac{1}{2}\leq \mu<1$. To this end, we need to estimate the term $\frac{A^\ast(y-y^\delta)}{\alpha}$ in (\ref{eq:tikh_sourcerep2}). It is, due to the noise model (\ref{eq:noise}),
\[
\frac{\|A^\ast(y-y^\delta)\|}{\alpha}\leq \frac{\|A^\ast\|\|y-y^\delta\|}{\alpha}\leq \frac{\delta}{\alpha}
\]
where we used $\|A\|=\|A^\ast\|=1$. Therefore, as long as $\frac{\delta}{\alpha}\rightarrow 0$, the noise component in (\ref{eq:tikh_sourcerep2}) vanishes as $\delta\rightarrow 0$.  This is the case for the well-known a-priori parameter choice $\alpha\sim \delta^{\frac{2}{2\mu+1}}$ (\ref{eq:a_priori_asc}), since for $\mu>\frac{1}{2}$, $\frac{\delta}{\alpha}\sim \delta^{\frac{2\mu-1}{2\mu+1}}\rightarrow 0$ as $\delta\rightarrow 0$.


\section{Higher order Tikhonov regularization}\label{sec:higherTikh}
In this section we illustrate the principles described in the previous section by considering higher order Tikhonov regularization and showing that any fixed saturation level can be achieved by slightly adjusting Tikhonov regularization in order to enforce higher smoothness of the approximate solutions $\xad$. It is known that \textit{iterated Tikhonov regularization} \cite{Gfrerer,KingChill} is able to do that, but the discrepancy principle fails to yield such rates \cite{Groetsch}, while it works with small restriction in our version. In addition to this, our method is direct and requires only the solution of a linear system.

Let $\kappa\in \N_0$. Then we calculate regularized approximations to (\ref{eq:problem}) as solution of
\begin{equation}\label{eq:tikhkappa}
(\basa^{\kappa+1}+\alpha I)x=A^\ast\basa^\kappa  y^\delta,
\end{equation}
i.e.,
\[
\xad=\mathrm{argmin}_{x\in X} \|\basa^{\frac{\kappa}{2}}(Ay-y^\delta)\|^2+\alpha\|x\|^2.
\]
This is similar to fractional Tikhonov regularization, see, e.g., \cite{filter2015,ReiHo}, where the idea is to reduce the smoothness of the approximate solutions. From the view of approximation, this is not a wanted effect, because we may decrease the approximation smoothness below the smoothness of $\xd$, which would yield sup-optimal convergence rates.

The solution to (\ref{eq:tikhkappa}) may, analogously to (\ref{eq:tikh_sourcerep1}), be written as
\begin{equation}\label{eq:tikhhigh_sourcerep}
\xad=\basa^{\kappa+\frac{1}{2}}\frac{y^\delta-A\xad}{\alpha},
\end{equation}
or in analogy to (\ref{eq:tikh_sourcerep2}),
\begin{equation}\label{eq:tikhhigh_sourcerep_apriori}
\xad=\basa^{\kappa+1}\frac{\xd-\xad}{\alpha}+\basa^\kappa A^\ast(y^\delta-y).
\end{equation}

\begin{figure}
\includegraphics[width=\linewidth]{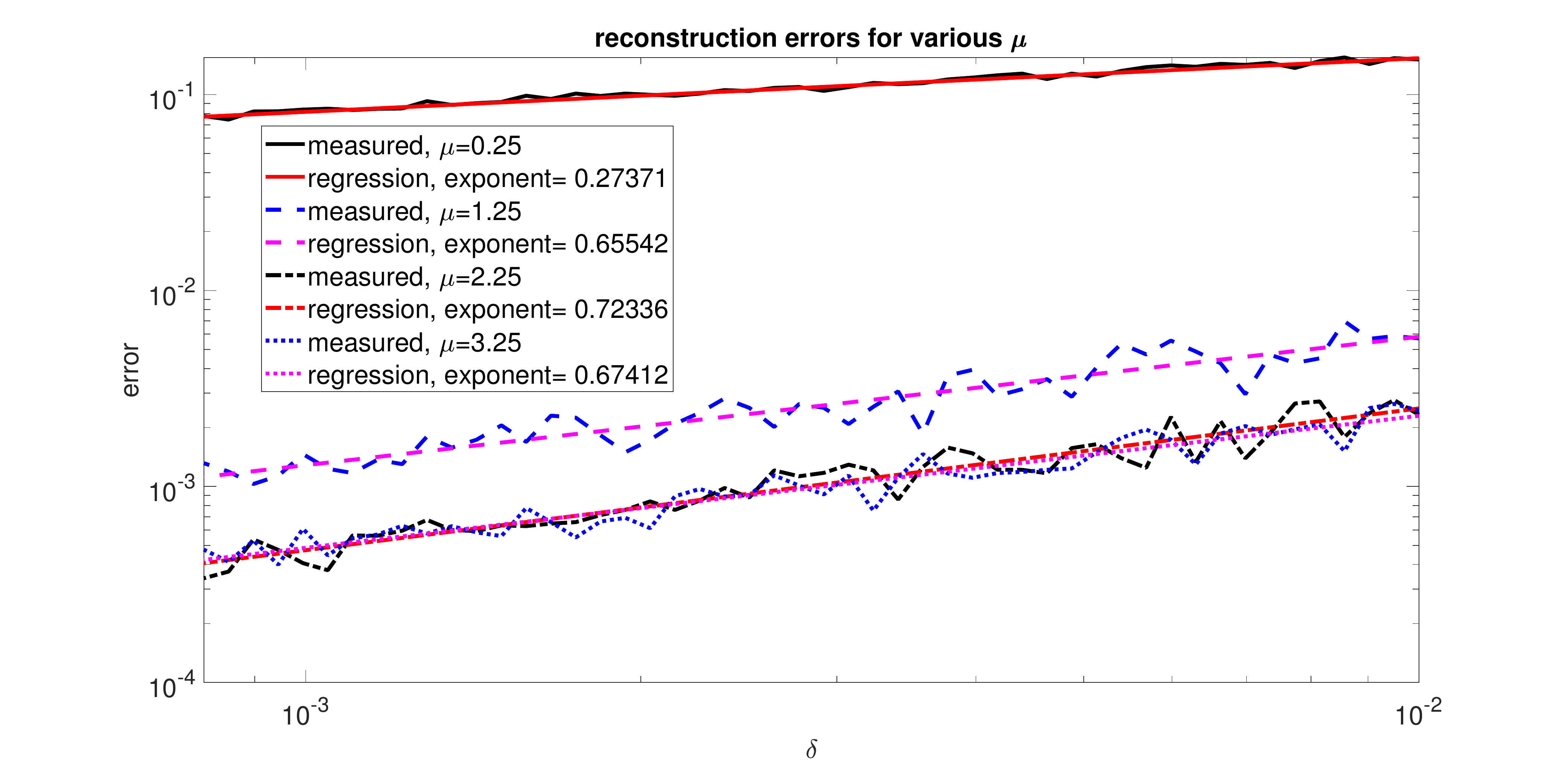}\caption{High-order Tikhonov regularization with $\kappa=1$, several solutions $\xd$ with differing $\mu$, and the discrepancy principle. For $\mu<1.5$ we obtain the order-optimal convergence rates with slight deviations to be expected numerically. For larger $\mu$, the observed convergence rates stagnates or even declines.}\label{fig:highorder}
\end{figure}

From (\ref{eq:tikhhigh_sourcerep}) we would read a saturation $\mu<\kappa+\frac{1}{2}$ for the discrepancy principle, since until then $\xad$ is smoother than $\xd$. As with classical Tikhonov regularization, we would expect an a-priori choice to saturate at $\mu=\kappa+1$ see (\ref{eq:tikhhigh_sourcerep_apriori}). Let as before $\xa$ denote the high-order Tikhonov approximation (\ref{eq:tikhkappa}) with noise-free data. With the singular system of $A$ we can write
\[
\xad=\sum_{\sigma_i>0} \frac{\sigma^{2\kappa+\frac{1}{2}}}{\sigma^{2\kappa+2}+\alpha}\langle y^\delta,u_n\rangle v_n.
\]
Now standard calculus of filter-based regularization, see. e.g. \cite{Louis}, yields
\begin{equation}\label{eq:error_hightikh}
\|\xad-\xd\|\leq \|\xa-\xd\|+\|\xad-\xa\|\leq c_1\alpha^{\frac{\mu}{\kappa+1}}+c_2\frac{\delta}{\alpha^{\frac{1}{2\kappa+2}}},
\end{equation}
where for the approximation error $\|\xa-\xd||$ one finds the condition $\mu<\kappa+1$, and for the noise amplification $\|\xad-\xa\|$ the requirement $\kappa>-\frac{1}{2}$ must hold. By balancing the terms in (\ref{eq:error_hightikh}) we obtain the a-priori parameter choice 
\[
\alpha=\tilde c \delta^{\frac{(\kappa+1)(2\kappa+2)}{(2\kappa+2\mu)+\kappa+1}}
\]
with an appropriate constant $\tilde c>0$, which yields the convergence rate
\[
\|\xad-\xd\|\leq c\delta^\frac{2\mu}{2\mu+1},\qquad \mu<\kappa+1
\]
just as expected. One can also show that indeed the discrepancy principle indeed saturates at $\mu=\kappa+\frac{1}{2}$, for example via \cite[Theorem 5.3.2]{Louis}. We demonstrate the saturation in Figure \ref{fig:highorder}. Using (\ref{eq:tikhkappa}) with $\kappa=1$, we calculate approximations to four solutions $\xd$ that satisfy a source condition (\ref{eq:sc_sum}) with $\mu\in \{ 0.25,1.25,2.25,3.25\}$. The regularization parameter is chosen according to the discrepancy principle, such that we expect a saturation of the convergence rate at $\mu=1.5$. The numerically observed exponents $q$ for the convergence rate prototype $\|\xad-\xd\|=\bigo(\delta^q)$ (in brackets the theoretical rate $\frac{2\mu}{2\mu+1}$) are $0.27 \, (0.33)$ for $\mu=0.25$, $0.65 \, (0.71)$ for $\mu=1.25$, $0.72 \, (0.81)$ for $\mu=2.25$, and $0.67 \, (0.86)$ for $\mu=3.25$.  For $\mu<1.5$ the numerical values are in line with the theoretical ones. After that, they no longer increase significantly with $\mu$.

As a final remark, we mention that such a high-order regularization reduces the influence of the noise, since the smoothing operator $A^\ast(A^\ast A)^\kappa$ is applied to $y-y^\delta$, where the smoothing becomes stronger the larger $\kappa$ is. The price to pay is that the approximation of $\xd$ in the smoother spaces will lead to the solutions looking too smooth such that the practical user will likely not be interested in unnecessarily large values of $\kappa$. From the regularization theoretical point of view, however, arbitrarily high $\kappa>0$ will yield order-optimal convergence rates.

\section{Oversmoothing regularization}\label{sec:os}
At this point we would like to comment on \textit{oversmoothing regularization}. The term appeared in recent years and is used for situations in which Tikhonov-type regularization functionals
\begin{equation}\label{eq:tikg_gen}
\|Ax-y^\delta\|_p^p+\alpha\Omega(x)
\end{equation}
with $1\leq p<\infty$ and a suitable penalty functional $\Omega:X\rightarrow \R$ are minimized under the specific assumption that the true solution $\xd$ does not attain a finite penalty value, $\Omega(\xd)=\infty$. Most of the works consider a Hilbert space setting, more precisely, Hilbert scales. We are only aware of few reports on this scenario in a Banach space setting, namely for $\ell^1$-regularization \cite{BerndBanach,GerHof19,HohMil}. The term oversmoothing regularization is due to the observation that the regularized solutions $\xad$, i.e., minimizers of (\ref{eq:tikg_gen}), fulfill $\Omega(\xad)<\infty$, and therefore are much smoother than $\xd$. One of the goals of this article was to show that the principle ``regularized solutions $\xad$ are smoother than exact solution $\xd$'' is already prevalent when $\Omega(\xd)<\infty$. It appears likely that one can show that this principle is necessary for order-optimal convergence rates, but this is left as an open problem. Instead, we show that our new approximation-based theory allows to understand the case of \textit{oversmoothing regularization} much easier. In the context of classical Tikhonov regularization (\ref{eq:tikh_classic}) with penalty $\Omega(x)=\|x\|_X^2$, oversmoothing regularization $\|\xd\|_X=\infty$ requires to determine convergence rates in a norm weaker than the $X$-norm, as by construction $\|\xad-\xd\|_X=\infty$. Therefore one would need to build up the convergence theory anew based on the new base norm, which would be too lengthy. Instead, we consider Tikhonov regularization in Hilbert scales, which is a well established concept in regularization theory.

Most of the Hilbert scale works on oversmoothing regularization have been centred on nonlinear forward operators \cite{GHH,HofMat17,HofMatapriori,HofPlat}, but, to be consistent with the previous parts of the paper, we will stay within the linear framework $A:X\rightarrow Y$. Let as before $X$ be our (base) Hilbert space, and let $T:X\rightarrow X$ be a densely defined, unbounded,  linear, and self-adjoint operator for which $\|Tx\|_X\geq C\|x\|_X$. Then $T$ generates a Hilbert scale $\{X_\nu\}_{\nu\in \R}$, and we have $X_\nu =\mathcal{D}(B^{\nu})=\cR(T^{-\nu})$. It is $X=X_0$ and we set $\|x\|_\nu:=\|T^\nu x\|_X$.
We consider again the noise-free data $y=Ax^\dag$ and noisy data $y^\delta$, $\|y-y^\delta\|\leq \delta$. The approximations to $\xd$ are obtained as minimizers of the Tikhonov-functional
\begin{equation}\label{eq:tikh_OS}
\|Ax-y^\delta\|^2+\alpha \|x\|_s^2
\end{equation}
where $s\geq 0$ is a chosen parameter that regulates the smoothness of the penalty. Such functionals have been analysed by Natterer in \cite{Natterer}, and he proved convergence rates under an a-priori parameter choice. The only ingredient needed is a stability-type estimate
\begin{equation}\label{eq:stab}
c_1\|x\|_{-a}\leq \|Ax\|\leq c_2\|x\|_{-a}
\end{equation}
for all $x\in X$, $0<c_1\leq c_2$ and some $a>0$. Assuming $\xd \in \cR(T^{-p})$, i.e., $\|\xd\|_p<\infty$, Natter shows that 
\begin{equation}\label{eq:results_natterer}
\|\xad-\xd\|\leq c \delta^{\frac{p}{a+p}} \quad\mbox{ if }\quad \alpha=c\delta^{\frac{2(s+a)}{a+p}},
\end{equation}
provided that $s\geq \frac{p-a}{2}$, which we write
\begin{equation}\label{eq:scond_natterer}
p\leq 2s+a.
\end{equation}
The remarkable observation is that there is no issue with having $s$ larger than $p$, i.e., enforcing a penalty smoothness that is (arbitrarily far!) above the smoothness of $\xd$. Only too low parameters $s$ are excluded via (\ref{eq:scond_natterer}). In other words there is no difference in the convergence rate and parameter choice between \textit{oversmoothing} regularization $\|\xd\|_s^2=\infty$ and ``traditional'' regularization $\|\xd\|_s^2<\infty$. The question is: Why? And how can one interpret this intuitively? Our new approach allows to answer these question. As before, the reason lies in the first order optimality condition. For (\ref{eq:tikh_OS}) this reads
\begin{equation}\label{eq:tikh_os_firstorder}
A^\ast(A\xad-y^\delta)+\alpha T^{2s}\xad=0,
\end{equation}
i.e.,
\begin{equation}\label{eq:tikhos_reordered}
\xad=T^{-2s}A^\ast \frac{A\xad-y^\delta}{\alpha}.
\end{equation}
The smoothness condition (\ref{eq:stab}) implies 
\begin{equation}\label{eq:ranges}
\cR(T^{-a})=\cR(A^\ast)=\cR(\basa^{\frac{1}{2}}),
\end{equation}
see \cite{BoetHof}. Therefore we can write $\xd\in \cR(T^{-p})=\cR(\basa^{\frac{p}{2a}})$ and, due to (\ref{eq:tikhos_reordered}),
\begin{equation}\label{eq:range_HS}
\xad\in \cR(T^{-(2s+a)})=\cR((A^\ast A)^{\frac{s}{a}+\frac{1}{2}}).
\end{equation}
We find again the smoothing principle we encountered previously for classical Tikhonov regularization. It is not only the penalty functional that defines the solution smoothness, it is its interaction with the forward operator, more precisely its adjoint. The actual solution smoothness is higher (when $\cR(A)\neq \overline{\cR(A)}$) than implied purely by the penalty. Namely, the penalty $\|\cdot\|_s^2$ yields $\|\xad\|_s<\infty$. However (\ref{eq:range_HS}) implies $\|\xad\|_{2s+a}<\infty$. The solutions are more than twice as smooth as implied by the penalty. The condition (\ref{eq:scond_natterer}) now states that optimal convergence rates are obtained whenever $\xad$ is at least as smooth as $\xd$, but $\xad$ can have any higher smoothness. Again, this is the same principle as with classical Tikhonov regularization discussed in the previous sections. In this sense, there is nothing special about \textit{oversmoothing regularization}. In fact, due to (\ref{eq:range_HS}) and (\ref{eq:ranges}), we can write the Tikhonov functional (\ref{eq:tikh_OS}) as
\begin{equation*}\label{eq:tikh_os_2}
\|\basa^{\frac{s}{a}}(Ax-y^\delta)\|^2+\alpha \|x\|^2
\end{equation*}
which yields a first-order optimality condition equivalent to (\ref{eq:tikh_os_firstorder}) but is not \textit{oversmoothing} whenever $p>0$. 

We can now recover Natterers results (\ref{eq:results_natterer}) through the approximation based approach. The main ingredient is Theorem \ref{thm:conv_approx} which we apply with $\nu\in\{0,\frac{1}{2}\}$, $\kappa=\frac{s}{a}+\frac{1}{2}$ and $\mu=\frac{p}{2a}$. Keeping these parameter assignments fixed, the results of Section \ref{sec:appsc} still hold and yield precisely (\ref{eq:results_natterer}). As a side result we obtain, practically for free, that the discrepancy principle also yields the rate from (\ref{eq:results_natterer}), and that furthermore it coincides up to constants with the a-priori choice as long as (\ref{eq:scond_natterer}) holds.

\section{Connection to Landweber iteration}\label{sec:lw}
A main message of this paper is that the source element in the approximation space is the driving element of regularization. To make this more tangible, we now demonstrate that Landweber iteration and Tikhonov regularization, an iterative regularization method and a variational regularization method, can, as long as $\xd$ is not too smooth such that Tikhonov regularization has saturated, be somewhat unified by this principle.

\begin{figure}
\includegraphics[width=0.49\linewidth]{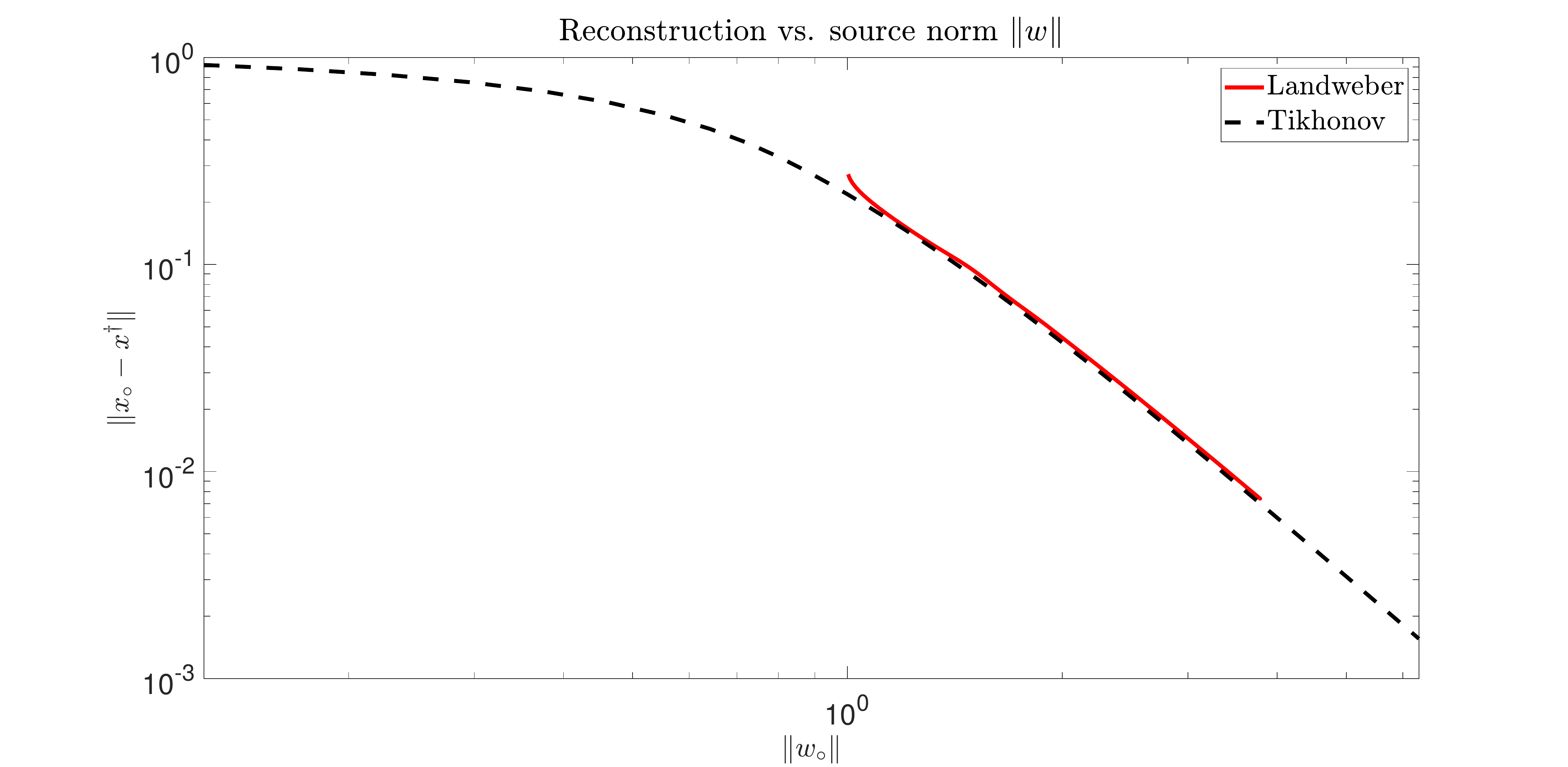}\includegraphics[width=0.49\linewidth]{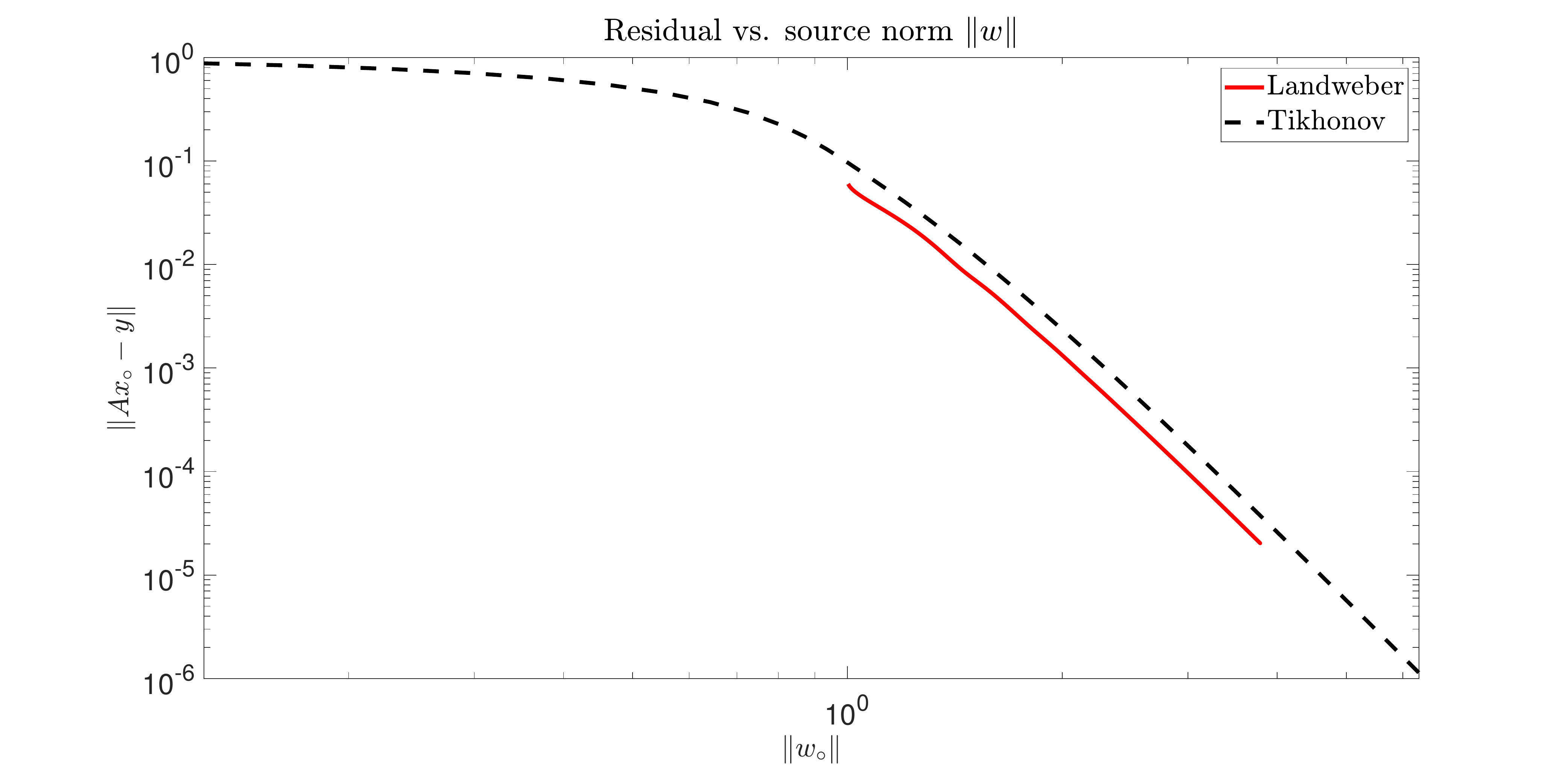}
\includegraphics[width=0.49\linewidth]{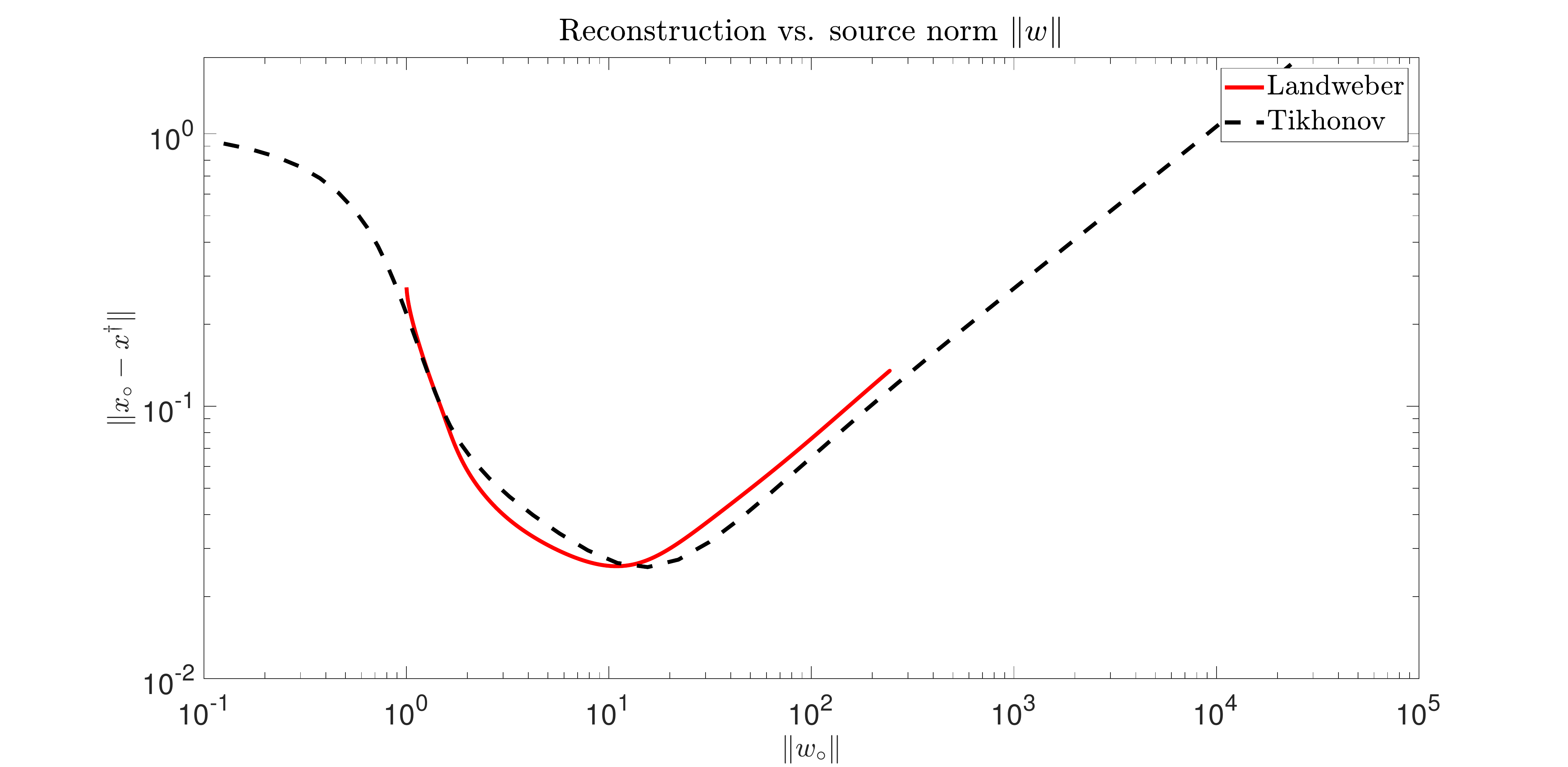}\includegraphics[width=0.49\linewidth]{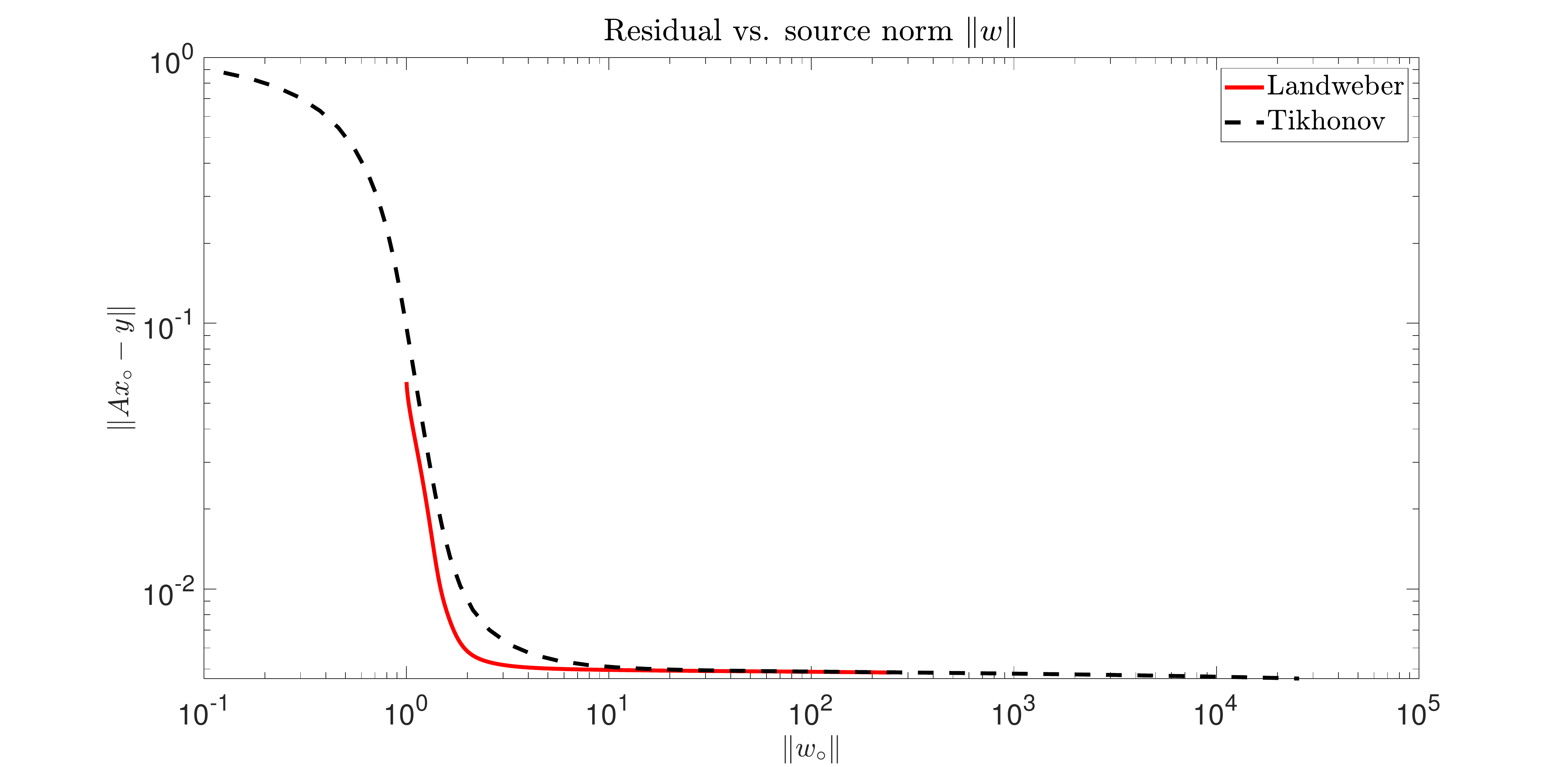}
\caption{Reconstruction error (left column) and residuals (right column) as functions of the source norm $\|w\circ\|$ for Tikhonv regularization (\ref{eq:tikh_x0} (black, dashed) and Landweber iteration (\ref{eq:lw} (red, solid). Top row: noise free, bottom row: $\delta=0.001$. Both methods yield almost identical graphs, indicating that the growth of the source element is crucial for the convergence analysis.}\label{fig:tikhlw}
\end{figure}

Landweber iteration is a classical iterative regularization method aimed at minimizing $\|Ax-y^\delta\|^2$ over $x\in X$. Starting from some point $x_0\in X$ one iterates
\[
x_{k+1}=x_k-\beta A^\ast(Ax_k-y^\delta)
\]
where the step-length parameter fulfils $0<\beta<\frac{2}{\|A\|^2}$ and $k=0,1,2,\dots$ until some chosen stopping index $k_\ast$. One can sum up the iteration from $x_0$ to $x_{k_\ast}$ to obtain
\begin{equation}\label{eq:lw}
x_{k_\ast}=x_0+ A^\ast \left(\beta\sum_{i=0}^{k_\ast} (I-\beta\asa)^i(y^\delta-Ax_0) \right)
\end{equation}
(see e.g. \cite[Eq. (6.3)]{EHN} with $\beta=1$ and $x_0=0$). As comparison, Tikhonov regularization (\ref{eq:tikh_classic}) with penalty $\|x-x_0\|^2$ yields the first-order condition
\begin{equation}\label{eq:tikh_x0}
x_{Tikh}=x_0+A^\ast\frac{A\xad-y^\delta}{\alpha}.
\end{equation}
Let in the following $\circ$ be a place holder for the indices $k_\ast$ and \textit{Tikh}. Comparing (\ref{eq:lw}) and (\ref{eq:tikh_x0}), we see that the difference between the methods is the way the source element $w_\circ$ in the source representation 
\begin{equation*}
x_\circ =x_0+A^\ast w_\circ
\end{equation*} 
is formed. Consequently, the regularization parameters $\alpha$ and $k_\ast$ serve as a way to control the growth of the source element, and in turn, the source element is crucial for the reconstruction error. To illustrate this, we conduct the following experiment: For a simple diagonal operator (\cite[Model Problem with $\eta=2$ and $\beta=2$ ($\xd$ satisfies (\ref{eq:sc_sum}) with $\mu=0.375$)]{GerRamres}) we calculate Tikhonov approximations $\xad$ and Landweber iterates $x_{k_\ast}$ for several parameters $\alpha$ and $k_{\ast}$ for noise free data $y$ and noisy data $y^\delta$ with $\delta=0.001$. We then plot the reconstruction error $\|x_\circ-\xd\|$ and the residual $\|Ax_\circ-y^\delta\|$ as functions of the norm of the source element $\|w_\circ\|$. The result, displayed in Figure \ref{fig:tikhlw}, is that the graphs are very close to each other and that they share the same convergence behaviour. We used 60 values $\alpha= 0.7^j$, $j=1,2,\dots,60$ and 50000 Landweber iterations. The graphs show that Tikhonov regularization moves through the approximation space $\cR(A^\ast)$ much more efficiently, as it covers a much wider range of source norms $\|w\|$.

Of course, this ``source similarity'' between Tikhonov regularization and Landweber iteration only holds if $\xd$ fulfils a source condition with $\mu<\frac{1}{2}$. If a source condition with $\mu>\frac{1}{2}$ is satisfied, the saturation for Tikhonov regularization set in. It is, however, well-known that Landweber iteration does not suffer from saturation, and that it yields optimal convergence rates for any $\mu>0$. We can explain this intuitively from the theory of source elements. Landweber iteration is a Krylov subspace method, i.e., the $k$-th iterate lies in the $k$-th Krylov subspace $\mathcal{K}_k(A^\ast y^\delta,\asa)$. It has a representation 
\[
x_k=x_0+A^\ast w_0^{(k)} +A^\ast \asa w_1^{(k)}+ A^\ast \basa^2 w_2^{(k)}+\dots+A^\ast \basa^k w_k^{(k)}
\]
with source source elements $w_i^{(k)}$, $i=0,1,\dots,k$ changing with each iteration. This means that for any fixed solution smoothness (\ref{eq:sc}), (\ref{eq:sc_open}) with $0<\mu<\infty$, the Landweber iterates are smoother than $\xd$ whenever $k>\mu$. As explained for Tikhonov regularization with $0<\mu<\frac{1}{2}$, this effect yields optimal convergence rates and a coincidence of a-priori parameter choice and discrepancy principle.

\section{Conclusions and outlook}
We have investigated classical Tikhonov regularization from a novel angle based on the smoothness properties of the approximate solutions, which can be interpreted as having a source condition for the regularized solutions. We demonstrated in this paper that this new approach can be used to understand regularization more intuitively. We recovered and explained the well-known convergence rate results and saturation. A main concept is that the regularized solutions have to be smoother than the one to be approximated. To demonstrate the simplicity of this principle, we proposed a variant of Tikhonov regularization can overcome any fixed saturation limit. We also showed that ``oversmoothing" regularization is easily explained through the approximation approach. Further, we showed that the idea of assigning a source condition to the regularized solutions allows to show that Landweber iteration and Tikhonov regularization function on the same basic principle, namely an appropriate choice of the norm of the source element.
This paper constitutes the first step in developing a more homogeneous theory for the regularization of ill-posed problems. The results are promising but more details have to be worked out. A major step is a generalization to nonlinear forward operators. In this situation, the derivative of the forward operator plays a prominent role in the first-order optimality condition, so that it seems likely to replace the role of $A$ in this paper with that derivative. Further, a generalization to Banach spaces is needed. One can then, in general, no longer use the classical source conditions. Instead, the first-order optimality conditions yield a smoothness condition for the subdifferential of the regularized solutions, which has then to be connected with $\xd$. It is an interesting question whether the principle of approximating $\xd$ through smoother objects still holds in this setting. Note that for the example of $\ell^1$-regularization this can be confirmed, since under mild assumptions the role the adjoint $A^\ast$ in the first-order condition ensures that the regularized solutions are finite dimensional, which is much stronger than just $\ell^1$-smoothness \cite{GerFle}.

\section*{Acknowledgements}
D. Gerth was supported by Deutsche Forschungsgemeinschaft (DFG), project GE3171/1-1 (Project Number 416552794) and would like to thank Prof. Oliver Ernst (TU Chemnitz) and Prof. Bernd Hofmann (TU Chemnitz) for the fruitful discussions.


\bibliographystyle{plain}

\end{document}